\documentclass[review]{siamltex}

\usepackage{amssymb,color,amsmath,bbm,latexsym,eucal,cite}
\usepackage{graphicx,psfrag} 
\usepackage{url}
\usepackage{epstopdf}

\newtheorem{ipotesi}{Assumption}[section]    

\newtheorem{remark}[theorem]{Remark}
\newtheorem{definizione}[theorem]{Definition}

\def\IR{\hbox{\rm I\kern-.2em\hbox{\rm R}}}
\def\IN{\hbox{\rm I\kern-.2em\hbox{\rm N}}}
\newcommand{\Pbk}{\mathbb{P}_{k-1}}

\def\t{\theta}
\def\Nk{N_k}
\def\Ntilde{\widetilde N_{k+1}}
\def\Nkp1{N_{k+1}}
\def\Dkp1{D_{k+1}}
\def\Nkpp{N_{k+1}^t}
\def\Nkppa{N_{k+1}}
\def\Dkp1{D_{k+1}}
\def\pred{{\rm{Pred}}_k}
\def\ared{{\rm{Ared}}_k}

\def\hit{{\cal{K}}_{\epsilon}}

\newcommand{\comment}[1]{}

\newcommand{\be}{\begin{equation}}
\newcommand{\ee}{\end{equation}}

\newcommand{\eqdef}{\stackrel{\rm def}{=}}

\newcounter{algo}[section] 
\renewcommand{\thealgo}{\thesection.\arabic{algo}}
\newcommand{\algo}[3]{\refstepcounter{algo}
\begin{center}\begin{figure}[!htbp]
\framebox[\textwidth]{
\parbox{0.95\textwidth} {\vspace{\topsep}
{\bf Algorithm \thealgo : #2}\label{#1}\\
\vspace*{-\topsep} \mbox{ }\\
{#3} \vspace{\topsep} }}
\end{figure}\end{center}}

\title{A stochastic first-order trust-region method
with inexact restoration for finite-sum minimization\footnotemark[3]
}
\author{ Stefania Bellavia\footnotemark[1], Nata$\check{{\rm s}}$a Kreji\'c\footnotemark[2], Benedetta Morini\footnotemark[1], Simone Rebegoldi\footnotemark[1]}

\begin{document}

\maketitle
\footnotetext[1]{Dipartimento  di Ingegneria Industriale, Universit\`a degli Studi di Firenze,
Viale G.B. Morgagni 40,  50134 Firenze,  Italia. Members of the INdAM Research Group GNCS. Emails:
stefania.bellavia@unifi.it, benedetta.morini@unifi.it, simone.rebegoldi@unifi.it}
\footnotetext[2]{ Department of Mathematics and Informatics, Faculty of Sciences, University of Novi Sad, Trg Dositeja
Obradovi\'ca 4, 21000 Novi Sad, Serbia, Email: natasak@uns.ac.rs. }
\footnotetext[3]{ 
The research that led to the present paper was partially supported by  
a grant of the group GNCS of INdAM and partially developed within the  Mobility Project: "Second order methods for optimization problems in Machine Learning" (ID:
RS19MO05) executive programme of Scientific and Technological cooperation between the Italian
Republic and the Republic of Serbia 2019-2022. The work of the second author was  supported by Serbian
Ministry of Education, Science and Technological Development, grant no. 451-03-9/2021-14/200125. The fourth author acknowledges financial support received by the IEA CNRS project VaMOS.}

\begin{abstract}
We propose a stochastic first-order trust-region method with inexact function and gradient evaluations 
for solving finite-sum minimization problems.  Using a suitable
reformulation of the given problem, our method combines the inexact restoration approach 
for constrained optimization with the trust-region procedure and random models. Differently from other recent stochastic trust-region schemes, our proposed algorithm improves feasibility and optimality in a modular way.  We provide the expected number of iterations for reaching  a near-stationary point by imposing some probability
accuracy requirements on random functions  and gradients  
which are, in general, less stringent than the corresponding ones in literature. We validate the proposed algorithm on some nonconvex optimization problems arising in binary classification and regression, showing that it performs well in terms of cost and accuracy, and allows to reduce the burdensome tuning of the hyper-parameters involved. 
\end{abstract}

{\small
\textbf{Keywords}: finite-sum minimization, inexact restoration, trust-region methods, subsampling, worst-case iteration complexity.
}

\begin{AMS}
65K05, 90C26, 68T05.
\end{AMS}


\section{Introduction}
In this paper we consider the finite-sum minimization problem 
\begin{equation}\label{minf}
\min_{x\in \IR^n}f_N(x)= \frac{1}{N}\sum_{i=1}^N \phi_i(x),
\end{equation}
where $N$ is very large and finite and $\phi_i: \IR^n\rightarrow \IR$, $1\le i\le N$, are continuously differentiable. 
A number of important problems can be stated in this form, e.g., classification problems in machine learning, data fitting problems, sample average approximations of an objective function given in the form of mathematical expectation. 
In recent years the need for efficient methods for solving (\ref{minf}) resulted in a large body of literature and a number of methods have been proposed and analyzed,  see e.g., 
the reviews  \cite{springer_survey, curtis_survey, nocedalsurvey}.  

It is common to employ subsampled approximations of the objective function and its derivatives  with the aim   of reducing the computational cost. Focusing on  first-order methods, the  stochastic gradient {\cite{RM}} and more contemporary variants like SVRG \cite{GSBR,jz},  SAG \cite{sag}, ADAM \cite{adam}  and SARAH \cite{sarah}  are widely used for their simplicity and   low cost per-iteration.  
They do not call for function  evaluations   but require tuning  the learning rate and further possible hyper-parameters such as the mini-batch size. Since the tuning effort may be very computationally demanding \cite{CS}, more sophisticated approaches 
use   stochastic   linesearch or trust-region strategies to adaptively choose the learning rate, see 
\cite{ BGM, springer_survey, BGMT, bcms, ChenMeniSche18, CS, PS}.  In this context, function and gradient approximations have to satisfy sufficient accuracy requirements with some probability. This, in turn, in case of approximations via sampling, requires adaptive choices of the sample sizes used.

In a further stream of works, problem (\ref{minf}) is reformulated as a constrained optimization problem and 
the sample size is computed deterministically using the Inexact Restoration (IR) approach.   
The IR approach has been successfully combined with either the linesearch strategy \cite{km} 
or the trust-region strategy \cite{bkm1,bkm2,bkm}; in these papers,
function and gradient estimates are built with gradually increasing accuracy and averaging on the same sample. 

We propose a novel trust-region method with random models based on the IR methodology. 
In our proposed method, feasibility and optimality are improved in a modular way, and the resulting procedure differs from the existing stochastic trust-region schemes \cite{BSV,bcms, BGMT2, ChenMeniSche18,XY}
in the acceptance rule for the step. We provide a theoretical analysis and give a bound
on the expected iteration complexity to satisfy  an approximate
first-order optimality condition; this calls for accuracy conditions on random gradients
that are assumed to hold with some sufficiently large but fixed probability and 
are, in general, less stringent than the corresponding ones in  \cite{BSV,bcms, BGMT2, ChenMeniSche18,XY}. 
Our theoretical analysis  improves over the  one for the stochastic trust-region method 
with inexact restoration given in \cite{bkm}, since we no longer rely on standard theory for deterministic unconstrained optimization
invoked  eventually when functions and gradients are computed exactly.

The paper is organized as follows.  In Section 2 we give an overview of random models employed 
in the trust-region framework and introduce the main features of our contribution. 
The new algorithm is proposed in Section 3 and studied theoretically with respect to the iteration complexity analysis. Extensive numerical results are presented in Section 4.

\section{Trust-region method with random models}\label{sec:2}
Variants of the standard trust-region method based on the use of random models have been 
presented, to our knowledge, in \cite{bkm, BSV,BGMT2, bcms, ChenMeniSche18,CSK,XY}. 
They consist in the adaptation of the trust-region framework to the case where random estimates of the derivatives are introduced and function values are either computed exactly \cite{BSV} or replaced by stochastic estimates \cite{bkm,bcms, BGMT2, ChenMeniSche18,CSK,XY}.

The computation and acceptance of the iterates parallel  the standard trust-region mechanism, and  the success of the procedure relies
on function values and models being sufficiently accurate with fixed and large enough probability.
The accuracy requests in the mentioned works show many similarities; here we illustrate some issues related to 
the works \cite{bcms,ChenMeniSche18, XY}, which are closer to our approach.

Let  $\|\cdot\|$ denote  the 2-norm throughout the paper.
At iteration $k$ of a first-order stochastic trust-region model,  
given $x_k$, the positive trust-region radius $\delta_k$ and a random approximation $g_k$ of $\nabla f_N(x_k)$, 
let consider the  model
$$
\varsigma_k(x_k+s)=f_N(x_k)+ g_k^T s $$
 for $f_N$ on $B(x_k, \delta_k)=\{x\in \mathbb{R}^n: \|x-x_k\|\le \delta_k\}$ and the trust-region  
problem $\min_{\|s\|\le \delta_k} \varsigma_k(x_k+s)$. Thus, the trust region step takes the form $s_k=-\delta_k g_k/\|g_k\|$.
 
Two estimates  $f^{k,0}$ and $f^{k,s}$ of $f_N$ at $x_k$ and  $x_k+s_k$, respectively, are employed to either accept or reject the trial point $x_k+s_k$.
The classical ratio between the actual and predicted reduction is replaced by
\begin{equation}\label{rho}
 \rho_k=\displaystyle \frac{f^{k,0}-f^{k,s}}{\varsigma_k(x_k)-\varsigma_k(x_k+s_k)},
 \end{equation}
and a successful iteration is declared when $\rho_k \ge \eta_1$  and $ \|g_k\|\ge \eta_2 \delta_k$ for some constants $ \eta_1 \in (0,1) $ and  positive and possibly large $ \eta_2 $. Note that the computation of both the step $s_k$ and  the denominator in (\ref{rho}) are independent of $f_N(x_k)$. 
Furthermore, note that a 
successful iteration might not yield an actual reduction in $f_N$ because the quantities involved in $\rho_k$ are random approximations to the true value of the objective function.

The condition $ \|g_k\|\ge \eta_2 \delta_k$ is not typical of standard trust-region and depends
on the fact that  $\delta_k$ controls the accuracy of function and gradients.
Specifically, the models used are required to be sufficiently accurate with some probability. 
The model $\varsigma_k$ is supposed to be, $p_M$-probabilistically, a $\kappa_*$-fully linear model of $f_N$ on the ball $B(x_k, \delta_k)$,   i.e., the requirement
\begin{equation}\label{model_acc}
|f_N(y)-\varsigma_k(y)| \le \kappa_* \delta_k^2, \quad \|\nabla f_N(y)-g_k\|\le \kappa_* \delta_k,  \quad y \in B(x_k, \delta_k) 
\end{equation}
with $\kappa_*>0$, has to be fulfilled at least with probability $p_M \in (0,1)$. 
Moreover, the estimates $f^{k,0}$ and $f^{k,s}$
are  supposed to be $p_f$-probabilistically  $\epsilon_F$-accurate estimates of $f_N(x_k)$
and  $f_N(x_k+s_k)$, i.e., the requirement
\begin{equation}\label{f_acc}
|f^{k,0} -f_N(x_k)|\le \epsilon_F \delta_k^2, \quad  |f^{k,s}-f_N(x_k+s_k)|\le  \epsilon_F \delta_k^2,
\end{equation}
has to be fulfilled  at least with probability $p_f\in (0,1)$. Clearly, if  $f_N$ is computed exactly then condition (\ref{f_acc}) is trivially satisfied.

Convergence analysis in \cite{bcms,ChenMeniSche18, XY} shows that for  $p_M $ and  $p_f $ sufficiently large it holds $\lim_{k\rightarrow \infty} \delta_k=0$  almost surely. 
Moreover, if $f_N$ is bounded from below and $\nabla f_N$
is Lipschitz continuous, then 
$\lim_{k\rightarrow \infty} \|\nabla f_N(x_k)\|=0$ almost surely. 
Interestingly, the accuracy in  (\ref{model_acc}) and (\ref{f_acc}) increases as the trust region radius gets smaller 
but the probabilities $p_M$ and $p_f$  are fixed.

For problem (\ref{minf}) it is straightforward to build approximations of   $f_N$ and  
$\nabla f_N$  by sample average approximations
\begin{equation}\label{subsample}
f_M(x)   =    \frac{1}{M}\sum_{i\in I_M}  \phi_i(x),\qquad \nabla f_S(x)   =    \frac{1}{S}\sum_{i\in I_S}  \nabla\phi_i(x),
\end{equation}
where $I_M$ and $I_S$ are subsets of   $\{1, \ldots, N\}$ of cardinality $| I_M|=M$  and $|I_S|=S$, respectively.
The choice of   sample size such that (\ref{model_acc}) and (\ref{f_acc}) hold in probability 
is discussed in \cite[\S 5]{ChenMeniSche18} as follows.
Let $\mathbb E[|\phi_i(x)-f_N(x)|^2]\le V_f$,  $\mathbb E[|\nabla \phi_i(x)-\nabla f_N(x)|^2]\le V_g$, $  i=1, \ldots, N$, with $\mathbb E$ being the expected value of a random variable, and assume  
\begin{equation}\label{samplesize_storm}
M\ge  \frac{V_f}{\epsilon_F^2(1-p_f)\delta_k^4},\;\;\;\;\; S \ge   \frac{V_g}{\kappa_*^2(1-p_g)\delta_k^2} \\
\mbox{ and } \max\{M,S\}\le N.
\end{equation}
Then  {  { $f^{k,0}$ and $f^{k,s}$ built as in (\ref{subsample}) with sample size $M$ satisfy  (\ref{f_acc}) with probability $p_f$, 
while $g_k$ built as in (\ref{subsample}) with sample size $S$ satisfies $\|\nabla f_N(x_k)-g_k\|\le \kappa_* \delta_k$
with probability $p_g$. Furthermore,
using Taylor expansion and Lipschitz continuity of $\nabla f_N$, it can be proved that
(\ref{model_acc}) is met with probability $p_M=p_f p_g$; consequently,
a  $\kappa_*$-fully linear model of $f_N$ in $B(x_k, \delta_k)$
is obtained. }
}

In principle, conditions (\ref{model_acc}), (\ref{f_acc}) and $\lim_{k\rightarrow \infty} \delta_k=0$ imply that 
$f^{k,0}$, $f^{k,s}$ and $g_k$ will be computed at full precision for $k$ sufficiently large.
On the other hand, 
in applications such as machine learning, reaching  full precision is unlikely since $N$ is very large and  termination is based
on the maximum allowed computational effort or on the validation error. 
\vskip 10pt

\subsection{Our contribution}
{ 
We propose a trust-region procedure with random models based on
(\ref{subsample}) and combine it with  the inexact restoration  (IR)
method for constrained optimization \cite{MP}. To this end,
we  make a simple transformation of (\ref{minf}) into a constrained problem. Specifically, 
letting $I_M$ be an arbitrary nonempty  subset of   $\{1, \ldots, N\}$ of cardinality $| I_M|$ equal to $M$, we 
reformulate problem (\ref{minf}) as 
\begin{equation} 
\begin{aligned}
 & \min_{x\in \IR^n}  f_M(x)   =    \frac{1}{M}\sum_{i\in I_M}  \phi_i(x),\\
 & \mbox{ s.t. }      M=N.    
\end{aligned}  \label{minf1}
\end{equation}

Using the IR strategy allows to improve feasibility and optimality in a modular way
and gives rise to a procedure that differs from the existing trust-region schemes
in the following respects. First, at each iteration a {\em reference} sample size is fixed and used as
a guess for the approximation of function values. Second, the acceptance rule for the step
{is based on} the condition $\|g_k\|\ge \eta_2 \delta_k$, for some $\eta_2>0$, and {a sufficient decrease condition on} a merit function that measures both  the reduction of the objective function
and the improvement in feasibility. Finally,  the expected iteration complexity to satisfy  an approximate
first-order optimality condition is given, provided that, at each iteration $k$, the gradient estimates satisfy 
accuracy requirements of order ${\cal{O}}\left(\delta_k\right)$; such accuracy requirements    implicitly govern   function approximations and   are, in general, less stringent than the corresponding ones in  \cite{BSV,bcms, BGMT2, ChenMeniSche18,XY}, as carefully detailed in Section \ref{sec:algo_analysis}.

Our theoretical analysis improves over the analysis carried out in \cite{bkm}
for a similar stochastic trust-region coupled with inexact restoration, since here
we do not rely on the occurrence of full precision, $M=N$ in (\ref{minf1}), reached eventually  and
do not  apply
standard theory for unconstrained optimization. In fact, the expected number of iterations until a prescribed accuracy is reached
is provided without invoking full precision.

}

\section{The Algorithm}\label{sec:algo_analysis}
{ In this section we introduce our new algorithm  referred to as SIRTR (Stochastic Inexact Restoration Trust Region).} 

{First, we introduce some issues of IR methods.} { The level of infeasibility 
with respect to the constraint $M=N$ in \eqref{minf1} is measured by the following function $h$.}
\begin{ipotesi}\label{assh}
Let $h:\{1,2,\ldots,N\}\rightarrow \IR$ be a monotonically decreasing function such that $h(1)>0$, $h(N)=0$.
\end{ipotesi}

This  assumption implies  that there exist some positive $\underline{h}$ and $\overline{h}$ such that
\begin{equation}\label{boundh}
\underline{h}\le h(M) \ \ \mbox{ if } \ \ 1\le M<N, \quad \mbox { and }\quad  h(M)   \le \overline{ h} \ \ \mbox{ if } \ \ 1\le M\le N.
\end{equation}
One possible choice is $h(M)=(N-M)/N, \   1\le M\le N$.

The IR methods improve feasibility and optimality in modular way using a merit function to balance the progress. Since the reductions in the objective function and infeasibility might be achieved to a different degree, the IR method employs the merit function 
\begin{equation}\label{merit}
 \Psi(x,M,\t)=\t f_M(x)+(1-\t)h(M),
\end{equation}
with $\t\in (0,1). $

Our SIRTR algorithm   is a trust-region method that employs first-order random models. 
At a generic iteration $k$, 
 we fix a {\it trial} sample size $\Nkpp$ and 
build a linear model $m_k(p)$  around $x_k$  of the form
\begin{equation}\label{model}
m_k(p)=f_{\Nkpp}(x_k)+g_k^Tp,
\end{equation}
{ where  $g_k$  is a random estimator to $\nabla f_N(x_k)$.}
Then, we consider the  trust-region  problem
\begin{eqnarray}\label{tr_pb}
\min_{ \|p\|\le \delta_k} m_k(p),
\end{eqnarray}
whose solution is 
\begin{equation} \label{direction}
p_k = - \delta_k \frac{g_k}{\|g_k\|}.
 \end{equation}
As in standard trust-region methods, we distinguish between successful  and unsuccessful iterations. However, we do not employ here the classical acceptance condition, but a more elaborate one that involves the merit function \eqref{merit}. 

{The proposed method is sketched in Algorithm \ref{IRTR_algo} and  its steps are now discussed.}
At a generic iteration $k$, we have at hand the outcome of the previous iteration:  {the iterate} $x_k$, the sample sizes 
$ \Nk $ and $\widetilde{N}_{k}$, the penalty parameter $ \t_k $, the flag  {\tt iflag}. If {\tt iflag=succ} the previous iteration was successful, i.e., $x_{k}=x_{k-1}+p_{k-1}$, if {\tt iflag=unsucc} the previous iteration was unsuccessful, 
i.e., $x_k=x_{k-1}$. 

The scheduling procedure for generating the trial sample size $\Nkpp$ consists of Steps 1 and 2 of SIRTR. 
At Step 1, we determine a reference sample size $\widetilde{N}_{k+1}\leq N$.  If  {\tt iflag=succ}, then
the infeasibility measure $ h $ is sufficiently decreased as stated in (\ref{feas}).
If {\tt iflag=unsucc}, $\widetilde{N}_{k+1}$ is left unchanged from the previous iteration, i.e., 
$\widetilde{N}_{k+1}=\widetilde{N}_k$. We remark that  \eqref{feas} trivially implies $\widetilde{N}_{k+1}=N$ if $N_k=N$ and that it 
holds at each iteration, even when it is not explicitly enforced at Step 1 (see forthcoming Lemma \ref{lemma_hntilde}). In principle $\widetilde{N}_{k+1}$ could be the trial sample size but we aim at giving more freedom to the sample size selection process.  Thus, at Step 2,  we choose a trial sample size $ \Nkpp$ complying with condition \eqref{new2}. On the one hand, such a condition allows the choice $ \Nkpp < \Ntilde $ in order to reduce the computational effort; on the other hand, the choice $ \Nkpp \geq \Ntilde $  is also possible in order to satisfy specific accuracy requirements that will be specified later. 
When $ \Nkpp < \Ntilde $, condition \eqref{new2} rules the largest possible distance between $\Nkpp$ and $\widetilde{N}_{k+1}$ in terms of $ \delta_k $; in case $ \Nkpp \ge  \Ntilde$, \eqref{new2} is trivially satisfied.

At Step 3 we {form the  linear random model \eqref{model} and compute its minimizer.} 
Specifically, 
we fix the cardinality $N_{k+1,g}$ and choose the set of indices   $I_{N_{k+1,g}}\subseteq\{1, \ldots, N\}$ of cardinality $N_{k+1,g}$.
Then, we compute  the estimator $g_k$  of $\nabla f_{N}(x_k)$ as 
\begin{equation}\label{gk}
g_k=
\frac{1}{N_{k+1,g}}\sum_{i\in I_{N_{k+1,g}}} \nabla \phi_i(x_k)
\end{equation}
and  the solution $p_k$ in (\ref{direction}) of the   trust-region subproblem (\ref{tr_pb}). Further, we {compute
 $m_k(p_k)$  where $m_k$ is defined in  \eqref{model} and }
\begin{equation}\label{f_algo}
f_{\Nkpp}(x_k)=\frac{1}{\Nkpp}\sum_{i\in I_{\Nkpp}} \phi_i(x_k),
\end{equation}
with  $I_{\Nkpp}\subseteq\{1, \ldots, N\}$ being a set of  cardinality $\Nkpp$.


{At Step 4 we compute the new penalty term $ \t_{k+1}. $  The computation relies on the predicted reduction defined as 
\begin{equation}\label{defpred}
\pred(\t)=\t(f_{\Nk}(x_k)-m_k(p_k))+(1-\t)(h(\Nk)-h(\Ntilde)),  
\end{equation}
where $\t\in(0,1)$. This predicted reduction is a convex combination of the usual predicted 
reduction $f_{\Nk}(x_k)-m_k(p_k) $  in trust-region methods,  and the predicted reduction  $ h(\Nk)-h(\Ntilde)  $ 
in infeasibility obtained in Step 1.  
The new parameter $ \t_{k+1} $ is computed  so that 
\begin{equation} \label{pred}
 \pred(\t) \geq \eta_1 (h(\Nk)-h(\Ntilde)).
 \end{equation}
If (\ref{pred}) is satisfied at  $\t=\t_k$} then $\t_{k+1}=\t_k$, otherwise $ \t_{k+1} $
is computed as the largest value for which the above inequality   holds (see forthcoming Lemma \ref{lemmatheta}).
 
{Step 5 establishes if the iteration is successful or not.
To this end, given a point $ \hat{x} $ and $\theta \in (0,1) $, 
the actual reduction of $\Psi$ at the point $ \hat{x}$ has the form
 \begin{eqnarray}\label{defared}
   \ared(\hat{x},\t) & = & \Psi(x_k, \Nk, \t )- \Psi(\hat{x}, \Nkpp, \t ) \nonumber  \\
  & = & \t (f_{\Nk}{(x_k)} -f_{\Nkpp}(\hat{x}))+(1-\t)(h(\Nk)-h(\Nkpp)), \qquad
\end{eqnarray} 
and the iteration is successful whenever the following two conditions are both satisfied
\begin{align}
\ared(x_k+p_k,\t_{k+1}) &\ge \eta_1 \pred(\t_{k+1})\label{eq:check}\\
 \|g_k\|&\ge \eta_2 \delta_k. \label{eq:additional}
\end{align}
Otherwise the iteration is declared unsuccessful.} 
If the iteration is successful, we accept the step and the trial sample size, set 
{\tt iflag=succ} and  possibly increase the trust-region radius through (\ref{Delta_update}); the
upper bound $ \delta_{\max} $ on the trust region size is imposed in (\ref{Delta_update}).
In case of unsuccessful iterations,  we reject both the step and the trial sample size, set 
{\tt iflag=unsucc} and decrease the trust region size. 

{Concerning conditions \eqref{eq:check} and \eqref{eq:additional}, 
we observe that the former mimics the classical acceptance criterion of standard trust-region methods  while the latter
drives $\delta_k$ to zero as $\|g_k\|$ tends to zero. 

\algo{IRTR_algo}{The  Stochastic IRTR algorithm }{
Given   $x_0\in \IR^n$, $N_0$ integer in $(0,N]$,  $\t_0\in (0,1)$,   
 $ 0 < \delta_0 < \delta_{\max}$, \\
$\gamma > 1,\,  r, \eta_1, \in (0,1),  \;   \mu, \, \eta_2>0$.
\vskip 4pt
\noindent
0. Set $k=0$, {\tt iflag}={\tt succ}. \\
1. If {\tt iflag}={\tt succ}\\
\hspace*{20pt}{Find} $\Ntilde$ such that $\Nk\leq \Ntilde\le N$ and 
\begin{equation}\label{feas}
h(\Ntilde)\le r h(\Nk),
\end{equation}
\phantom{1.} {Else set $\Ntilde= \widetilde{N}_k$}. \\
 2. If $N_k=N$ set $N_{k+1}^t=N$\\
 \phantom{1.} Else find $\Nkpp$ such that
\begin{eqnarray} \label{tilde}
 h(\Nkpp)-h(\Ntilde) &\le& \mu \delta_k^{2}.\label{new2}
\end{eqnarray}
3. Choose $N_{k+1,g}$,  $I_{N_{k+1,g}}\subseteq\{1, \ldots, N\}$ s.t. $|I_{N_{k+1,g}}|=N_{k+1,g}$.\\
\phantom{1.} Compute $g_k$ as in (\ref{gk}), and set 
$$
p_k =- \delta_k \frac{g_k}{\|g_k\|}.
$$
\phantom{1.} Compute $f_{N_{k+1}^t}(x_k)$ as in (\ref{f_algo}), and $m_k(p_k)=f_{N_{k+1}^t}(x_k)+g_k^Tp_k$.\\
4.  Compute the penalty parameter $ \theta_{k+1}$
\begin{equation}
\begin{aligned}
\t_{k+1}=
\left\{
\begin{array}{ll}
&\t_k  \hspace*{30pt} \mbox{  if } \  \pred(\t_k)\ge\eta_1( h(\Nk)-h(\Ntilde))\\
& \displaystyle \frac{(1-\eta_1)(h(\Nk)-h(\Ntilde))}{m_k(p_k)-f_{\Nk}(x_k)+h(\Nk)-h(\Ntilde)} \quad  \mbox{otherwise}.
\end{array}
\right.
\end{aligned}\label{tkp1}
\end{equation}
5.  If $\ared(x_k+p_k,\t_{k+1})\ge \eta_1 \pred(\t_{k+1})$ and $\|g_k\|\ge \eta_2 \delta_k$ ({\em successful iteration})\\
\hspace*{20 pt} define
\begin{align}
x_{k+1}&=x_k+p_k\nonumber\\
\delta_{k+1}&=\min \left\{ \gamma \delta_k, \delta_{\max}\right\}\label{Delta_update}
\end{align}
 \hspace*{20 pt}set  {$N_{k+1}=N_{k+1}^t$}, $k=k+1$, {\tt iflag=succ} and go to Step 1.\\
\hspace*{10pt}
Else ({\em unsuccessful iteration}) define
\begin{align}
x_{k+1}&=x_k\nonumber\\
\delta_{k+1}&= \frac{\delta_k}{\gamma}\label{Delta_update_u}
\end{align}
    \\
\hspace*{20 pt}  set {$\Nkppa=N_k$},    $ k=k+1$, {\tt iflag=unsucc} and go to Step 1.
\caption{}}

We conclude the description of Algorithm \ref{IRTR_algo} showing that condition \eqref{feas} holds for all iterations, even when it is not explicitly enforced at Step 1.

\begin{lemma}\label{lemma_hntilde}
 Let Assumption 3.1 holds and $r\in (0,1)$ be the scalar in Algorithm \ref{IRTR_algo}.
The sample sizes $\Ntilde\leq N$ and $\Nk\leq N$ generated by  Algorithm \ref{IRTR_algo} satisfy
\begin{equation}\label{feas_bis}
h(\Ntilde)\le r h(\Nk), \quad \forall k\ge 0 .
\end{equation}
\end{lemma}
{\em Proof.} 
We observe that, by Assumption 3.1, (\ref{feas_bis}) trivially holds whenever $\Nk=\Ntilde=N$.   

Otherwise, we proceed by induction. Indeed, the thesis trivially holds for $k=0$, 
as we set {\tt iflag=succ} at the first iteration and  enforce  \eqref{feas_bis} at Step 1. 
Now consider a generic iteration $\bar k\ge 1$ and suppose
that  \eqref{feas_bis} holds for $\bar k-1$. If iteration $\bar k-1$ is successful, then condition \eqref{feas_bis} is enforced for iteration $\bar{k}$ at Step 1. 

If iteration $\bar{k}-1$ is unsuccessful, then at Step 5 we  set $N_{\bar k}=N_{\bar k-1}$. Successively,
at Step 1 of iteration $\bar{k}$ we set  $\widetilde N_{\bar k+1}=\widetilde N_{\bar k}$. 
Since \eqref{feas_bis} holds by induction at iteration $\bar k-1$, we have $h(\widetilde N_{\bar k})\le r h(N_{\bar k-1})$, which can be rewritten as $h(\widetilde N_{\bar k+1})\le r h(N_{\bar k})$ due to the previous assignments at Step 5 and Step 1. Then condition \eqref{feas_bis} holds also at iteration $\bar k$.  
$\Box$

\subsection{On the sequences    $\{\t_k\}$  and $  \{\delta_k\}$}
In this section, we analyze the properties of Algorithm \ref{IRTR_algo}. 
In particular, we prove that the sequence $\{\t_k\}$ is non increasing and 
uniformly bounded from below, and that 
the trust region radius $\delta_k$ tends to $0$ as $k\rightarrow \infty$.
We make the following assumption.
\vskip 5pt
\begin{ipotesi}\label{assxk}
Functions $\phi_i$ are continuously differentiable for $i=1,\ldots,n$.
There exist $\Omega\subset \mathbb{R}^n$ and $f_{low}$, $f_{up}$ such that 
$$
f_{low}<f_{M}(x)<f_{up}, \quad   \ 1\leq M\leq N, \ x\in\Omega, 
$$
and all iterates generated by Algorithm \ref{IRTR_algo} belong to $\Omega$.
\end{ipotesi}
\vskip 5 pt
 In the following, we let
\begin{equation}\label{kappaphi}
\kappa_\phi=   \max\{ |f_{low}|, |f_{up}|\}.
\end{equation}
\begin{remark}\label{remark_ml}
In the context of machine learning, the above assumption is verified in several cases, e.g., the mean-squares loss function coupled with either the sigmoid,  the softmax or the hyperbolic tangent  activation function; the mean-squares loss function coupled with ReLU or ELU  activation functions and proper bounds on the unknowns; the logistic loss function coupled with proper bounds on the unknowns \cite{Good}.
\end{remark}
\vskip 5pt
In the analysis that follows we will consider two options for $ \hat{x}$  in \eqref{defared}, $\hat{x} = x_k + p_k $  for successful iterations and $ \hat{x} = x_k $ for unsuccessful iterations.

Our first result characterizes the  sequence $\{\t_k\}$ of the penalty parameters; the proof follows closely \cite[Lemma 2.2]{bkm}.
\begin{lemma}\label{lemmatheta} 
Let Assumptions \ref{assh} and \ref{assxk} hold. 
Then the sequence $ \{\theta_k\} $  is positive, non increasing and bounded from below, $\t_{k+1}\ge \underline{\t}>0$ with $\underline{\t}$ independent of $k$ and 
\eqref{pred} holds with $\theta=\theta_{k+1}$. 

\end{lemma}

{\em Proof.} {We note that $\theta_0>0$ and proceed by induction assuming that $\theta_k$ is positive.
Due to Lemma \ref{lemma_hntilde}, for all iterations $k$ we have that $N_k\leq \widetilde{N}_{k+1}$ and 
that $N_k=\widetilde{N}_{k+1}$ if and only if $N_k=N$.}
First consider the case where $\Nk=\Ntilde$ (or equivalently $\Nk=\Ntilde=N$); then it holds
$h(\Nk)-h(\Ntilde)=0$, and $\Nkpp=N$ by Step 2.
Therefore, we have $\pred(\t)=\t \delta_k \|g_k\|>0$ for any positive  $\t$, and (\ref{tkp1}) implies  $\t_{k+1}=\t_k $.

Let us now  consider the case $\Nk<\Ntilde$. If inequality $\pred(\t_k) \ge  \eta_1(h(\Nk)-h(\Ntilde))$ holds  then 
(\ref{tkp1}) gives $\theta_{k+1}=\theta_k$.  Otherwise,  we have 
$$
\t_k  \left( f_{\Nk}(x_k)-m_k(p_k)-( h(\Nk)-h(\Ntilde) )\right)   < {(\eta_1-1)\left(h(\Nk)-h(\Ntilde)\right)} ,
$$
and since the right hand-side   is negative by assumption, it follows
$$
f_{\Nk}(x_k)-m_k(p_k)-(h(\Nk)-h(\Ntilde))<0.
$$ 
Consequently, $\pred(\t )\ge \eta_1(h(\Nk)-h(\Ntilde))$ is satisfied if 
$$
\t (f_{\Nk}(x_k)-m_k(p_k)-(h(\Nk)-h(\Ntilde)))\ge  (\eta_1-1)(h(\Nk)-h(\Ntilde)),
$$
i.e., if 
$$ \t \le \t_{k+1}\eqdef \frac{(1-\eta_1)(h(\Nk)-h(\Ntilde))}{m_k(p_k)-f_{\Nk}(x_k)+h(\Nk)-h(\Ntilde)}. $$
Hence $\t_{k+1}$ is the largest value satisfying (\ref{pred}) and $ \t_{k+1} < \t_k. $

Let us now prove that $ \t_{k+1} \geq \underline{\t}. $   
Note that by (\ref{feas_bis}) and (\ref{boundh})  
\begin{equation}\label{riduh}
h(\Nk)-h(\Ntilde)\ge (1-r)h(\Nk)\ge (1-r) \underline{h}.
\end{equation} 
Using (\ref{kappaphi})  
\begin{eqnarray*}
 m_k(p_k)-f_{\Nk}(x_k)+h(\Nk)-h(\Ntilde) &\le& m_k(p_k)-f_{\Nk}(x_k)+h(\Nk)\\
& \le&   f_{\Nkpp}(x_k)-\delta_k \|g_k\| - f_{\Nk}(x_k)+ \overline{h}  \\ 
 & \le&  |f_{\Nkpp}(x_k)-f_{N_k}(x_k)| +  \overline{h}    \le  2 k_{\phi}+ \overline{h},
\end{eqnarray*}
and  $\t_{k+1}$ in (\ref{tkp1}) satisfies 
\begin{equation}\label{eq:theta_sub}
 \t_{k+1} \ge \underline{\t}=\frac{(1-\eta_1)(1-r) \underline{h}}{2 k_{\phi}+ \overline{h}},
\end{equation}
which completes the proof.
$\Box$ 
\vskip 5pt
}

In the following, we derive bounds for the actual reduction $\ared(x_{k+1},\t_{k+1})$
in case of successful iterations and distinguish the iteration indexes $k$ as below:
\begin{eqnarray}
	{\cal I}_1&=&\{ k\ge 0 \mbox{  s.t.  }    N_k<\Ntilde \}  \label{set1},\\
	{\cal I}_2&=&\{ k\ge 0 \mbox{  s.t.  }  N_k=\Ntilde \} \label{set2}.
\end{eqnarray}
Note that 
$\mathcal{I}_1,\mathcal{I}_2$ are disjoint and any iteration index $k$ belongs to exactly one of these subsets. Moreover, 
 \eqref{feas_bis} yields $\Ntilde=N_k= \Nkpp=N$ for any $k \in \mathcal{I}_2$.

\begin{lemma} \label{case123}
Let Assumptions \ref{assh}-\ref{assxk} hold and suppose that iteration $k$ is successful. If   $k \in {\cal I}_1$ then
{
\begin{equation}\label{succ3}
\ared(x_{k+1},\t_{k+1})   \geq 
\frac{\eta_1^2(1-r)\underline{h}}{\delta_{\max}^2}\delta_k^2.
\end{equation}}
Otherwise, 
\begin{equation}\label{succ4}
\ared(x_{k+1},\t_{k+1})  \geq \eta_1\eta_2 \underline \theta\delta_k^2.
\end{equation}
\end{lemma}
{\em Proof.} Since iteration $k$ is successful,  $x_{k+1}=x_k+p_k$ and \eqref{eq:check} hold. Suppose $k \in {\cal I}_1$. By  \eqref{eq:check}  and \eqref{pred} 
$$ \ared( x_k + p_k ,\t_{k+1}) \geq \eta_1 \pred(\t_{k+1}) \geq \eta_1^2 (h(N_k) - h(\widetilde{N}_{k+1})). $$  
In virtue of Lemma \ref{lemma_hntilde} we have $ h(N_k) - h(\Ntilde) \geq (1-r) h(\Nk)$, hence we obtain
$$ \ared(x_k+p_k, \t_{k+1}) \geq \eta_1^2 (1-r) h(\Nk). $$ 
Dividing and multiplying the right-hand side above by $\delta_{k}^2$, applying the inequalities $ \underline{h}\leq h(N_k)$, $\delta_k\leq \delta_{\max}$, we get (\ref{succ3}).

Suppose $k \in {\cal I}_2$.  Then $N_k=\Ntilde$ and by  the definition of
$\pred(\t_{k+1})$  and Lemma \ref{lemmatheta}, we have 
$$
\pred(\t_{k+1}) = \t_{k+1}(f_N(x_k) - m_k(p_k)) = \t_{k+1} \delta_k \|g_k\| \geq \underline{\t} \delta_k\|g_k\|,
$$
and  therefore  \eqref{eq:check}, (\ref{eq:additional}) and Lemma \ref{lemmatheta}  yield (\ref{succ4}).
$ \Box$ 
\vskip 5pt

Let us now define a Lyapunov type function $ \Phi $   inspired by the paper  \cite{ChenMeniSche18}.
Assumption \ref{assh} implies that $ h(N_k) $ is bounded from above while Assumption \ref{assxk} implies that 
$ f_{N_k}(x) $ is bounded from below if $x\in \Omega$. Thus, 
there exists a constant $ \Sigma$ such that 
\begin{equation}\label{sigma}
f_{N_k}(x) - h(N_k) + \Sigma \geq 0 , \quad  x  \in \Omega, \quad k\ge 0.
\end{equation}
\vskip 5pt
\begin{definizione}\label{def_Phi}
Let  $ v \in (0,1) $ be a fixed constant. 
For all $k\geq 0$, we define
 \be \label{lf}
 \phi_k\eqdef\Phi(x_k,N_k,\t_k, \delta_k) = v\left(\Psi(x_k,N_k,\t_k) + \t_k \Sigma\right ) + (1-v) \delta_k^2,
 \ee
 where $\Psi$ is the merit function given in \eqref{merit} and $\Sigma$ is given in (\ref{sigma}).
\end{definizione}
\vskip 5pt
The choice of $v\in (0,1)$ in the above definition will be specified below. 
First, note that $\phi_k$ is bounded below for all $k\ge 0$, 
\begin{eqnarray} \label{philow}
\phi_k&\ge&   v\left(\Psi(x_k,N_k,\t_k) + \t_k \Sigma\right )  \nonumber \\
&\ge& v\left(\t_k f_{N_k}(x_k)+(1-\t_k)h(N_k)+\t_k (-f_{N_k}(x_k)+h(N_k))\right )\nonumber\\
&\ge& v   h(N_k) \ge 0.
\end{eqnarray}

Second,  adding and subtracting suitable terms, by the definition (\ref{lf}) and for all $k\geq0$,  we have 
 \begin{eqnarray}
 \phi_{k+1} - \phi_{k} & = &   v\left( \t_{k+1}f_{N_{k+1}}(x_{k+1}) + (1-\t_{k+1})h(N_{k+1})\right) \nonumber \\  
 &    & -v\left( \t_{k}f_{N_{k}}(x_{k}) + (1-\t_{k})h(N_{k})\right)  + v(\t_{k+1} - \t_k)\Sigma + (1-v)(\delta_{k+1}^2 - \delta_k^2) \nonumber\\
& = &  v\left( \t_{k+1}f_{N_{k+1}}(x_{k+1}) + (1-\t_{k+1})h(N_{k+1})\right) \pm v \t_{k+1} f_{N_{k}}(x_{k}) \pm v(1-\t_{k+1}) h(N_k) \nonumber\\
 &    & -v\left( \t_{k}f_{N_{k}}(x_{k}) + (1-\t_{k})h(N_{k})\right)  + v(\t_{k+1} - \t_k)\Sigma + (1-v)(\delta_{k+1}^2 - \delta_k^2) \nonumber \\
&=& v\left({\theta_{k+1}(f_{N_{k+1}}(x_{k+1})-f_{N_k}(x_k))+(1-\theta_{k+1})(h(N_{k+1})-h(N_k))}\right) \nonumber\\
 &    & + v(\t_{k+1} - \t_k)(f_{N_k}(x_k) - h(N_k) +\Sigma) +(1-v)(\delta_{k+1}^2 - \delta_k^2). \label{phiared}
\end{eqnarray}

{If the iteration $k$ is successful, then using \eqref{sigma}, the monotonicity of $\{\theta_k\}_{k\in\mathbb{N}}$ 
proved in Lemma \ref{lemmatheta}, and the fact that $N_{k+1}=N_{k+1}^t$, the equality \eqref{phiared} yields
\begin{equation}\label{phiared1}
\phi_{k+1} - \phi_{k}\leq -v \ared(x_{k+1},\theta_{k+1})+(1-v)(\delta_{k+1}^2 - \delta_k^2).
\end{equation} 
Otherwise, if the iteration $k$ is unsuccessful, then $x_{k+1}=x_k$, $N_{k+1}=N_k$ and thus the first quantity at the right-hand side of equality \eqref{phiared} is zero. Hence using again \eqref{sigma} and the monotonicity of $\{\theta_k\}_{k\in\mathbb{N}}$, we obtain
\begin{equation}\label{phiared1_unsuccess}
\phi_{k+1} - \phi_{k}\leq (1-v)(\delta_{k+1}^2 - \delta_k^2).
\end{equation}
} 
  
Now we provide  bounds for the change of $\Phi$ along subsequent iterations and again
distinguish the two cases $k\in {\cal I}_1, {\cal I}_2$ stated in (\ref{set1})-(\ref{set2}). 
\vskip 5pt
\begin{lemma}
Let Assumptions \ref{assh}-\ref{assxk} hold.
\begin{description}
\item{i)} {If the iteration $k$ is unsuccessful, then}
\be \label{phiared_BU}
 {\phi_{k+1} - \phi_k  \leq \chi_1 \delta_k^2, \qquad \chi_1= (1-v)\frac{1-\gamma^2}{\gamma^2}.} 
\ee
\item{ii)}  {If the iteration $k$ is successful and $k\in\mathcal{I}_1$, then}
\be \label{phiared_AS}
 {\phi_{k+1} - \phi_k  \leq \chi_2 \delta_k^2, \qquad \chi_2=
 \left({-v\left(\frac{\eta_1^2(1-r)\underline{h}}{\delta_{\max}^2}\right)}+(1-v)(\gamma^2-1)  \right).}
\ee
\hspace*{-20pt}
If  the iteration $k$ is successful and $k\in\mathcal{I}_2$, then 	
\be \label{phiared_BS}
 {\phi_{k+1} - \phi_k  \leq \chi_3\delta_k^2, \qquad \chi_3= \left(- v \eta_1\eta_2 \underline \theta  + (1-v)(\gamma^2-1)  \right) .}
\ee
\end{description}
\end{lemma}

{\em Proof.}  

{i)  If iteration $k$ is unsuccessful, the updating rule \eqref{Delta_update_u} for $\delta_{k+1}$
 implies $\delta_{k+1}= \delta_k/\gamma$. Thus, equation \eqref{phiared1_unsuccess} directly yields \eqref{phiared_BU}.  

ii)  If iteration $k$ is successful, the updating rule \eqref{Delta_update} for $\delta_{k+1}$ implies $\delta_{k+1}\le \gamma \delta_k$.
Thus combining \eqref{phiared1} with Lemma \ref{case123} we obtain \eqref{phiared_AS} and  \eqref{phiared_BS}. 
$\Box$ 	
	
\vskip 5pt

We are now ready to prove that a sufficient decrease condition holds for $\Phi$ along subsequent iterations
and that $\delta_k$ tends to zero. 

\begin{theorem} \label{theorem_sigma}
Let Assumptions \ref{assh}--\ref{assxk} hold.
There exists  $\sigma>0$,  depending on $v\in(0,1)$ in (\ref{lf}), such that
\begin{equation}\label{ineq_sigma}
\phi_{k+1}- \phi_k\le -\sigma \delta_k^2, \quad  \mbox{ for all } k\ge 0.
\end{equation}
\end{theorem}

{\em Proof. } In case of unsuccessful iterations, \eqref{phiared_BU} provides a sufficient decrease
$\phi_{k+1}- \phi_k$ for any value of $v\in(0,1)$.

In case of  successful iterations, $\chi_2$ and $\chi_3$  in \eqref{phiared_AS} and \eqref{phiared_BS} are both negative if 
\begin{equation}\label{ineqv_2}
\max\left\{ {\frac{(\gamma^2-1)\delta_{\max}^2}{\eta_1^2(1-r)\underline{h}+(\gamma^2-1)\delta_{\max}^2}},\frac{\gamma^2-1}{\eta_1\eta_2 \underline \theta+\gamma^2-1}\right\}<v<1. 
\end{equation}
Therefore, if $v$ is chosen as above and 
\begin{equation}\label{sigma_choice}
\sigma=\min\{\chi_1,\, \chi_2,\, \chi_3\},
\end{equation}
then  \eqref{phiared_BU}--\eqref{phiared_BS} imply (\ref{ineq_sigma}) and the proof is completed.
$\Box$
\vskip 5pt
\begin{theorem}\label{theorem_delta}
Let Assumptions  \ref{assh}-\ref{assxk} hold. 
Then the sequence 
$\{\delta_k\}$ in Algorithm  \ref{IRTR_algo}  satisfies
$$
\lim_{k\rightarrow \infty} \delta_k=0.
$$
\end{theorem}
{\em Proof.} 
Under the stated conditions Theorem \ref{theorem_sigma} holds and   summing up \eqref{ineq_sigma} for $j=0,1,\ldots,k-1$, we obtain
\begin{equation*}\label{sumdelta}
  \phi_{k}-\phi_0= \sum_{j=0}^{k-1} ( \phi_{j+1}-\phi_j)\le -\sigma  \sum_{j=0}^{k-1}  \delta_j^2.
\end{equation*}
Given that, by \eqref{philow}, $ \phi_k $ is bounded from below  for all $ k, $   we conclude that 
$
\sum_{j=0}^{\infty}  \delta_j^2<\infty
$, 
and hence 
$
\lim_{j\rightarrow \infty}  \delta_j= 0.$
\subsection{Complexity analysis}

{Algorithm \ref{IRTR_algo} generates a random process since 
the function estimates $f_{N_{k+1}^t}(x_k)$ in \eqref{f_algo} and gradient estimates $g_k$ in \eqref{gk} are random.} 
All random quantities are denoted by capital letters, while the use of small letters is reserved for their realizations. In particular, the iterates $X_k$, the trust region radius $\Delta_k$,  the gradient estimates $G_k,\nabla f_{N_{k+1}^t}(X_k)$, and  
the value $\Phi_k$ of   the function $\Phi$ in (\ref{lf})  at iteration $k$
are random variables,  while $x_k$, $\delta_k$, $g_k$ and  $\phi_k$  are their realizations.  
We denote with $\mathbb{P}_{k-1}(\cdot)$ and $\mathbb{E}_{k-1}(\cdot)$ the probability and expected value conditioned to the past  until iteration $k-1$.

In this section, our aim is to derive a bound on the expected number of iterations   that occur in Algorithm \ref{IRTR_algo} to reach
a desired accuracy. We show that our algorithm is included 
into the stochastic framework given in \cite[\S 2]{bcms} and consequently we derive an upper bound on the expected value of 
the hitting time $\hit$ defined below.

\vskip 5pt
\begin{definizione}\label{Nepsilon}
Given $\epsilon>0$,  the hitting time $\hit$ is the random variable
\begin{equation*}
	\hit=\min\{k\geq 0: \ \|\nabla f_N(X_k)\|\leq \epsilon\},
\end{equation*}
i.e., $\hit$ is the first iteration such that $\|\nabla f_N(X_k)\|\le \epsilon$.
\end{definizione}
\vskip 5pt \noindent

Our analysis relies on the assumption that $g_k$ and $\nabla f_{\Nkpp}(x_k)$ are 
probabilistically accurate 
estimators of the true gradient at $x_k$, in the sense that 
the events
\begin{eqnarray}
{\cal{G}}_{k,1} &=&\{ \|\nabla f_N(X_k)-G_k\| \leq \nu \Delta_k  \}, \label{event_accu1}\\
{\cal{G}}_{k,2} &=&\{\|\nabla f_N(X_k) - \nabla f_{\Nkpp}(X_k)\| \leq \nu \Delta_k \}, \label{event_accu2}
\end{eqnarray}
are true at least with probability $\pi_1\in (0,1) $ and $\pi_2\in (0,1)$, respectively.  
Using the same terminology of \cite{BCS,Cartis-et-al}, we say that iteration $k$ is {\em true} if both $\mathcal{G}_{k,1}$ and $\mathcal{G}_{k,2}$ are true. Furthermore, we introduce the two random variables
\begin{equation}\label{eq:Ik}
	I_k=\mathbbm{1}_{\mathcal{G}_{k,1}}, \quad J_k=\mathbbm{1}_{\mathcal{G}_{k,2}},
\end{equation}
where $\mathbbm{1}_A$ denotes the indicator function of an event $A$.

Finally, we need the following additional assumptions.
\begin{ipotesi}\label{assf}
The gradients  
$ \nabla \phi_i$  are  Lipschitz continuous with constant $L_i$. Let ${L=\frac{1}{2}\max_{1\le i\le N} L_i }$.
\end{ipotesi} 
\vskip 2pt
\begin{ipotesi}\label{assg}
There exists $g_{\max}$ such that 
$$
\|g_k\|\le g_{\max}, \quad  k\ge 0 .
$$
\end{ipotesi} 
 We observe that the loss functions mentioned in Remark \ref{remark_ml} satisfy Assumption \ref{assg}. 

First, we  analyze the occurrence of successful iterations  and show that 
the availability of  accurate  gradients  has an impact on the acceptance of the trial steps.
The following lemma establishes  that 
if the iteration $k$  is {\em true} and $\delta_k$ is smaller than a certain threshold, then the iteration is successful.
The analysis is  presented for a single realization of Algorithm  \ref{IRTR_algo} 
and specializes for $k$ in the sets $ {\cal I}_1$, $ {\cal I}_2$.
\vskip 5pt
\begin{lemma} \label{succaccurate}
Let Assumptions \ref{assh}-\ref{assg} hold and suppose that iteration $k$ is true. 
\begin{description}
\item{i)} If $k \in {\cal I}_1$, then the iteration is successful whenever  
\begin{equation} \label{Delta_ASaccu}
\delta_k \le  \min \left\{  \frac{\|g_k\|}{\eta_3} , \, \frac{\|g_k\|}{\eta_2} \right\},
\end{equation}
where  $\eta_3=\frac{\delta_{\max}g_{\max}(\theta_0(2\nu+L)+(1-\underline{\theta})\mu)}{\eta_1(1-\eta_1)(1-r)\underline{h}}$. 
\item{ii)} If $k \in {\cal I}_2$, then the iteration  is successful whenever  
\be \label{Delta_ASaccuI2}
{ \delta_k \le \min \left\{ \frac{ (1-\eta_1)\|g_k\|}{2\nu+L}  ,\, \frac{\|g_k\|}{\eta_2} \right\}.}
\ee
\end{description}
\end{lemma}

{\em Proof. } From Assumption \ref{assf},  it follows that $\nabla f_{N_{k+1}^t}$ is Lipschitz continuous with constant $2L$. 
Then,  
\begin{align}
 |m_k(p_k) - f_{ {N_{k+1}^t}}(x_k + p_k)| & = \left |\int_{0}^{1} \left(g_k \pm \nabla f_{ {N_{k+1}^t}}(x_k) 
-\nabla  f_{ {N_{k+1}^t}}(x_k + \tau p_k) \right)^T p_k d \tau  \right | \nonumber \\
& \leq  \int_{0}^{1} \| g_k-\nabla  f_{{N_{k+1}^t}}(x_k)\| \|p_k\| d \tau  + \int_{0}^{1} 2L  \tau \|p_k\|^2  d \tau\nonumber\\
& \leq  \int_{0}^{1} ({\| g_k-\nabla  f_{N}(x_k)\|+\| \nabla f_N(x_k)-\nabla  f_{{N_{k+1}^t}}(x_k)\|}) \|p_k\| d \tau \nonumber\\
&+ \int_{0}^{1} 2L  \tau \|p_k\|^2  d \tau
\label{difffm}
\end{align} 
and, since ${\cal{G}}_{k,1}$ and ${\cal{G}}_{k,1}$ are both true, (\ref{event_accu1}) and (\ref{event_accu2}) yield 
\begin{equation}\label{modelestimate2}
|m_k(p_k) - f_{{N_{k+1}^t}}(x_k + p_k)|\le ({2}\nu+L )\delta_k^2.
\end{equation}
Now, let us analyze condition (\ref{eq:check})  for successful iterations.

i) If  $k \in {\cal I}_1$, by  (\ref{defpred}), (\ref{defared}) and (\ref{pred})  we obtain 
\begin{eqnarray} \label{aredcase1}
\ared(x_k+p_k,\t_{k+1}) - \eta_1 \pred(\t_{k+1}) &=& (1-\eta_1)\pred(\t_{k+1}) +\ared(\t_{k+1})-\pred(\t_{k+1}) \nonumber \\
& =&  (1-\eta_1)  \pred(\t_{k+1}) + \t_{k+1}(m_k(p_k)-f_{\Nkpp}(x_{k}+p_k)) \nonumber \\
&  +& (1-\t_{k+1})(h(\Ntilde)-h(\Nkpp)) \nonumber \\ 
 &  \geq & \eta_1 (1-\eta_1) (h(\Nk)-h(\Ntilde))\nonumber\\  
 &+&\t_{k+1}(m_k(p_k)  -f_{\Nkpp}(x_{k}+p_k)) \nonumber\\ 
& +&(1-\t_{k+1})(h(\Ntilde)-h(\Nkpp)). 
\end{eqnarray} 
Using \eqref{modelestimate2}, \eqref{new2} and {${\underline{\theta}}\leq \theta_{k+1}\leq \theta_0$}, we also have 
\begin{align} \label{bound1}
{\t_{k+1}(f_{\Nkpp}(x_{k}+p_k)-m_k(p_k)) }&{ +   (1-\t_{k+1})(h(\Nkpp)-h(\Ntilde))} \nonumber\\
 &  \leq (\theta_0(2\nu+L)+(1-\underline{\theta})\mu) \delta_k^2 . \end{align} 
{Note that the combination of \eqref{feas_bis}, \eqref{boundh}, \eqref{Delta_update}  and Assumption \ref{assg}, guarantees that
\begin{equation}\label{bound2}
h(N_k)-h(\Ntilde)\geq (1-r)h(N_k)\geq \frac{(1-r)\underline{h}\delta_k\|g_k\|}{\delta_{\max}g_{\max}}.
\end{equation} 
Then, from  (\ref{aredcase1}), (\ref{bound1}), and (\ref{bound2}), we have 
\begin{align*}
\ared(x_k+p_k,\t_{k+1}) - \eta_1 \pred(\t_{k+1})&\ge \frac{\eta_1(1-\eta_1)(1-r)\underline{h}\delta_k\|g_k\|}{\delta_{\max}g_{\max}} \\
&- { (\theta_0(2\nu+L)+(1-\underline{\theta})\mu) \delta_k^2}.
\end{align*} 
Combining this result with 
(\ref{eq:additional}), the proof is complete.}
\vskip 5pt

\vskip 5pt
ii) 
Using (\ref{defpred}), (\ref{defared}), $k \in {\cal I}_2$, we have
\begin{eqnarray*}
\ared(x_k+p_k,\t_{k+1}) - \eta_1 \pred(\t_{k+1})&=&(1-\eta_1)\pred(\t_{k+1}) +\ared(\t_{k+1})-\pred(\t_{k+1})\\
& = & (1-\eta_1) \t_{k+1}\delta_k\|g_k\|+\t_{k+1}(m_k(p_k)-{f_N}(x_k+p_k))
\end{eqnarray*}
Using (\ref{modelestimate2}) we get 
\begin{eqnarray} \label{eq:ared_pred_I2_I3}
\ared(x_k+p_k,\t_{k+1}) - \eta_1 \pred(\t_{k+1})  &\ge & (1-\eta_1)\t_{k+1} \delta_k\|g_k\| \nonumber \\ 
& -&\ \ \t_{k+1}(2\nu+L)\delta_k^2.
\end{eqnarray}
Combining the above inequality with 
(\ref{eq:additional}),  we have proved that the iteration is successful whenever   \eqref{Delta_ASaccuI2} holds.
$ \Box$ 
\vskip 5pt
We can now guarantee that a successful iteration $k$ occurs whenever $k$ is true, 
the prefixed accuracy $\epsilon$ in Definition  \ref{Nepsilon} has not been achieved at $k$, 
and $\delta_k$ is below a certain threshold depending on $\epsilon$. 
Again, the result is stated for a single realization of the algorithm.
{
\begin{lemma}\label{ddagger}
Let Assumptions \ref{assh}-\ref{assg} hold. Suppose that  $\|\nabla f_N(x_k)\|> \epsilon$, for some $\epsilon>0$, 
the iteration $k$ is true, and
\begin{equation} \label{def_b}
\delta_k<\delta^\dagger:= \min \left \{ \frac{\epsilon}{2\nu}, \frac{\epsilon}{2\eta_2},  \frac{\epsilon}{2\eta_3} ,\frac{ \epsilon(1-\eta_1)}{2(2\nu+L)} \right\}.
\end{equation}
Then, iteration k is successful.
\end{lemma}}
\vskip 3pt\noindent
\begin{proof}
By $\|\nabla f_N(x_k)\|>\epsilon$, the occurrence of ${\mathcal{G}_{k,1}}$   and (\ref{def_b}),  we have  
$$
\|g_{k}-\nabla f_N(x_k)\|\le \nu \delta_k <\frac{\epsilon}{2},
$$
and this yields $\|g_k\|\ge \frac{\epsilon}{2}$.  Then,  Lemma \ref{succaccurate} implies that 
iteration $k$ is successful.
\end{proof}
\vskip 5pt

We now proceed similarly to \cite[\S 2]{bcms} and analyse the random process $\{(\Phi_k,\Delta_k,W_k)\}_{k\in\mathbb{N}}$ generated by Algorithm \ref{IRTR_algo}, where $\Phi_k$ is  the random variable whose realization is given in (\ref{lf}) and $W_k$ is the random variable defined as
\begin{equation}\label{eq:Wk}
\begin{cases}
W_0 =1\\
W_{k+1}=2\left(I_kJ_k-\frac{1}{2}\right), \quad k=0,1,\ldots
\end{cases}
\end{equation}
Clearly, $W_k$ takes values $\pm 1$. Then, we can prove the following result.
\vskip 5pt
\begin{lemma}\label{ass_bcs}
Let Assumptions \ref{assh}-\ref{assg} hold, $v$ as in (\ref{ineqv_2}),  $\delta^\dagger$ as in (\ref{def_b})
and  $\hit$ as in Definition \ref{Nepsilon}. Suppose there exists some $j_{\max}\geq 0$ such that $\delta_{\max}=\gamma^{j_{\max}}\delta_0$, and $\delta_0>\delta^{\dagger}$. Assume that the estimators $G_k$ and $\nabla f_{N_{k+1}^t}(X_k)$ are independent random variables, and the events ${\cal{G}}_{k,1},{\cal{G}}_{k,2}$ occur with sufficiently high probability, i.e., 
	\begin{eqnarray}\label{eq:Gk-prob}
		\Pbk({\cal G}_{k,1})= \pi_1, \quad \Pbk({\cal G}_{k,2})= \pi_2, \quad \text{and }p=\pi_1\pi_2>\frac{1}{2}.
	\end{eqnarray}
	Then,
\begin{itemize}
\item[i)] there exists $\lambda>0$ such that $\Delta_k\leq \delta_0e^{\lambda \cdot j_{\max}}$ for all $k\geq 0$;
\item[ii)] there exists a constant $ \delta_{\epsilon}=\delta_0e^{\lambda \cdot j_{\epsilon}}$ for some $j_{\epsilon}\leq 0$ such that, for all $k\geq 0$,
\begin{equation}\label{eq:ii}
\mathbbm{1}_{\{\hit>k\}}\Delta_{k+1}\geq \mathbbm{1}_{\{\hit>k\}}\min\{\Delta_ke^{\lambda W_{k+1}}, \delta_{\epsilon} \},
\end{equation}
where  $W_{k+1}$  satisfies
\begin{equation}\label{eq:Wk-prob}
	\Pbk(W_{k+1}=1)=p, \quad \Pbk(W_{k+1}=-1)=1-p;
\end{equation}
\item[iii)] there exists a nondecreasing function $\ell:[0,\infty)\rightarrow (0,\infty)$ and a constant $\Theta>0$ such that, for all $k\geq 0$,
\begin{equation}
	\mathbbm{1}_{\{\hit>k\}}\mathbb{E}_{k-1}[\Phi_{k+1}]\leq \mathbbm{1}_{\{\hit>k\}}\Phi_k-\mathbbm{1}_{\{\hit>k\}}\Theta \ell(\Delta_k).
\end{equation}
\end{itemize}
\end{lemma}
{\em Proof.}
The proof parallels that of \cite[Lemma 7]{bcms}.	
	
i) Since $\delta_{\max}=\gamma^{j_{\max}}\delta_0$, we can set $\lambda=\log(\gamma)>0$, and the thesis follows from Step 5 of Algorithm \ref{IRTR_algo}.

ii) Let us set
\begin{equation}\label{eq:delta_epsilon}
	{ \delta_\epsilon=\frac{\epsilon}{\xi}}, \quad \text{where }\xi\geq \max\left\{2\nu,2\eta_2,
	{  2\eta_3,}\frac{2(2\nu+L)}{1-\eta_1}\right\},
\end{equation}
and assume that $ \delta_{\epsilon}=\gamma^{j_{\epsilon}}\delta_0$, for some integer $j_{\epsilon}\leq 0$; notice that we can always choose $\xi$ sufficiently large so that this is true. 
As a consequence, $\Delta_k= \gamma^{i_k} \delta_\epsilon$ for some integer $i_k$.

When $\mathbbm{1}_{\{\hit>k\}}=0$, inequality \eqref{eq:ii} trivially holds. Otherwise, conditioning on $\mathbbm{1}_{\{\hit>k\}}=1$, we can prove that
\begin{equation}\label{eq:thesis}
\Delta_{k+1}\geq \min\{\delta_\epsilon,\min\{\delta_{\max},\gamma\Delta_k\}I_kJ_k+\gamma^{-1}\Delta_k(1-I_kJ_k)\}.
\end{equation}
Indeed, 
for any realization such that $\delta_k>\delta_{\epsilon}$, we have $\delta_k\geq \gamma\delta_{\epsilon}$ and because of Step 5, it follows that $\delta_{k+1}\geq \delta_{\epsilon}$. Now let us consider a realization such that  $\delta_k\leq \delta_{\epsilon}$. Since $\hit>k$ and $\delta_{\epsilon}\leq \delta^{\dagger}$, if $I_kJ_k=1$ (i.e., $k$ is true), then we can apply Lemma \ref{ddagger} and conclude that $k$ is successful. Hence, by Step 5, we have $\delta_{k+1}=\min\{\delta_{\max},\gamma\delta_k\}$. If $I_kJ_k=0$, then we cannot guarantee that $k$ is successful; however, again using Step 5, we can write $\delta_{k+1}\geq \gamma^{-1}\delta_k$. Combining these two cases, we get \eqref{eq:thesis}. If we observe that $\delta_{\max}=\gamma^{j_{\max}}\delta_0\geq { \gamma^{j_{\epsilon}}\delta_{0}= \delta_{\epsilon}}$, and recall the definition of $W_k$ in \eqref{eq:Wk}, then equation \eqref{eq:thesis} easily yields \eqref{eq:ii}. 
The probabilistic conditions \eqref{eq:Wk-prob} are a consequence of \eqref{eq:Gk-prob}. 

(iii) The thesis trivially follows from \eqref{ineq_sigma} with $\ell(\Delta)=\Delta^2$ and $\Theta=\sigma$. $\Box$
\vskip 5pt

The previous lemma shows that the random process $\{(\Phi_k,\Delta_k,W_k)\}_{k\in\mathbb{N}}$ complies with Assumption 2.1 of \cite{bcms}.

\begin{theorem}\label{teo_E}
	Under the assumptions of Lemma \ref{ass_bcs}, we have
	\begin{equation}{ 
		\mathbb{E}[\hit]\leq \frac{p}{2p-1}\cdot \frac{\phi_0\xi^2}{{\sigma}  \epsilon^2}}+1.
	\end{equation}
	{where $\xi$ is chosen as in \eqref{eq:delta_epsilon}} { and $\sigma$ is given in \eqref{sigma_choice}.}
	\end{theorem}
	
{\em Proof.}
The claim follows directly by \cite[Theorem 2]{bcms}. $\Box$

\vskip 5 pt

\begin{remark}
The requirement of (\ref{event_accu1}) and (\ref{event_accu2}) to 
hold in probability is less stringent than the overall conditions (\ref{model_acc}) 
and (\ref{f_acc}).
Analogously to the discussion in Section \ref{sec:2}, 
if $\mathbb E[|\nabla \phi_i(x)-\nabla f_N(x)|^2]\le V_g$, $  
i=1, \ldots, N$, then Chebyshev inequality guarantees that 
events (\ref{event_accu1}) and (\ref{event_accu2}) hold in probability 
when
\begin{equation*}
\frac{V_g}{\nu^2(1-\pi_1)\delta_k^2} \leq N_{k+1,g} \leq N   , \quad \frac{V_g}{\nu^2(1-\pi_2)\delta_k^2 }\leq N_{k+1}^t \leq N.
\end{equation*}
Clearly, $\min\{N_{k+1,g}, N_{k+1}^t\}={\cal{O}}(\delta_k^{-2})$ and in general these sample sizes are 
expected to growth slower than in (\ref{samplesize_storm}).



Finally, the complexity theory presented improves on \cite{bkm} where the iteration complexity before reaching  full precision  
$M=N$ in (\ref{minf1}) is estimated, and  thereafter existing  iteration complexity results for trust-region methods 
applied to (\ref{minf}) are invoked.
\end{remark}

\section{Numerical experience}
In this section, we evaluate the numerical performance of SIRTR on  some  nonconvex optimization problems arising in binary classification { and regression}.

All the numerical results have been obtained by running MATLAB R2019a on an Intel Core i7-4510U CPU 2.00-2.60 GHz with an 8 GB RAM. For all our tests, we { equip SIRTR with} $\delta_0=1$ as the 
initial trust-region radius, $\delta_{\max}=100$,  $\gamma=2$, $\eta=10^{-1}$, $\eta_2=10^{-6}$. Concerning   the inexact restoration phase, we  borrow the implementation details from \cite{bkm}.
Specifically, the infeasibility measure $h$ and the initial penalty parameter $\theta_0$ are set as follows: 
$$
h(M)=\frac{N-M}{N}, \quad  \t_0=0.9.
$$
The updating rule for choosing $\Ntilde$ has the form
\begin{equation}\label{eq:ntilde}
\Ntilde=\min\{N,\lceil \widetilde{c} N_k\rceil\},
\end{equation}
where $1<\widetilde{c}<2$ is a prefixed constant factor; note that this choice of 
$\Ntilde$ satisfies \eqref{feas} with $r=(N-(\widetilde{c}-1))/N$. At Step 2 the function sample size $\Nkpp$ is computed using the rule 
\begin{equation}\label{eq:nk1}
{\Nkpp}=\begin{cases}
\lceil \Ntilde-\mu N \delta_k^2\rceil, \quad & \text{if }\lceil \Ntilde-\mu N\Delta_{k}^2\rceil\in[N_0,0.95N]\\
\Ntilde, \quad & \text{if }\lceil \Ntilde-\mu N\Delta_{k}^2\rceil< N_0\\
N, \quad & \text{if }\lceil \Ntilde-\mu N\Delta_{k}^2\rceil> 0.95 N.
\end{cases}
\end{equation}

Once the set $I_{\Nkpp}$  is fixed, the search direction $g_k\in\mathbb{R}^n$ 
is computed via sampling as in \eqref{gk} and the sample size 
$N_{k+1,g}$ is fixed as
\begin{equation}\label{eq:nkg}
N_{k+1,g}=\lceil c  \Nkpp \rceil,
\end{equation}
with $c\in(0,1]$ and $I_{N_{k+1,g}}\subseteq  I_{\Nkpp}$.

\subsection{SIRTR performance}

In the following, we show the numerical behaviour of SIRTR on {nonconvex binary classification problems}. Let $\{(a_i, b_i)\}_{i=1}^N$ denote the pairs forming a training set with $a_i \in \IR^n$ containing the entries of the $i$-th example, and $b_i\in\{0, 1\}$ representing the corresponding label. Then, we address the following minimization problem  
\begin{equation}\label{eq:sigmoid}
\min_{x\in\IR^n}f_N(x)=\frac 1 N\sum_{i=1}^N \left(b_i-\frac{1}{1+e^{-a_i^Tx}} \right)^2,
\end{equation}
where the nonconvex objective function $f_N$ is obtained by composing a least-squares loss with the sigmoid function.

In Table \ref{test}, we report the information related to the datasets employed, including the number $N$ of training examples, the dimension $n$ of each example and the dimension $N_T$ of the testing set $I_{N_T}$. 

\begin{table}[t!] 
\begin{center}
\begin{tabular}{l  |r r |r }
&     \multicolumn{2}{r|} { Training set }& \multicolumn{1}{r} { Testing set }\\
\hline 
Data set & $N$ & $n$  &  $N_T$   \\  \hline
{\sc A8a}\cite{UCI} & 15887& 123& 6809\\
{\sc A9a}\cite{UCI} & 22793& 123& 9768\\
{\sc Cina0}   &  10000 & 132 & 6033 \\
{\sc cod-rna} \cite{libsvm}  & 41675& 8& 17860  \\
 {\sc Covertype} \cite{UCI} & 464810& 54 & 116202 \\
 {\sc Htru2} \cite{UCI} & 10000  & 8  & 7898  \\
{\sc Ijcnn1}\cite{libsvm} & 49990 & 22& 91701\\
{\sc Mnist} \cite{mnist}  & 60000& 784& 10000  \\
{\sc phishing} \cite{libsvm}  & 7739& 68& 3316  \\
{\sc real-sim} \cite{libsvm}  & 50616 & 20958 & 21693  \\
{\sc w7a}\cite{libsvm} & 17284 & 300 & 7408\\
{\sc w8a}\cite{libsvm} & 34824& 300& 14925\\
\hline
\end{tabular}
\caption{Data sets used}\label{test}
\end{center}
\end{table} 

We focus on three aspects: the classification error provided by the final iterate, the computational cost, the occurrence of termination before full accuracy in function evaluations is reached. The last issue is crucial because it indicates the ability of the inexact restoration approach to solve (\ref{eq:sigmoid}) with random models and to rule
sampling  and steplength selection.

The average classification error provided by the final iterate, say $x_{{\rm fin}}$, is defined as  
\begin{equation}\label{class-err}
{\tt err}=\frac{1}{N_T}\sum\limits_{i\in I_{N_T}}|b_i-b_i^{pred}|,
\end{equation}
where $b_i$ is the exact label of the $i-$th instance of the testing set, and $b_i^{pred}$ is the corresponding predicted label, given by $b_i^{pred}=\operatorname{max}\{\operatorname{sign}(a_i^Tx_{{\rm fin}}),0\}$. 

The computational cost  is measured in terms of full function and gradient evaluations. 
In our test problems, the main cost in the computation of $\phi_i$, $1\le i\le N$,  is the scalar product $a_i^Tx$: once this product is evaluated, it can be reused for computing $\nabla \phi_i$. Nonetheless, following \cite[Section 3.3]{XuRoosMaho_proc}, we count both function and gradient evaluations as if we were addressing a classification problem based on a neural net. Thus, computing a single function $\phi_i$ requires $\frac{1}{N}$ forward propagations, whereas the gradient evaluation corresponds to $\frac{2}{N}$ propagations (an additional backward propagation is needed). { Note that, once $\phi_i$ is computed, the corresponding gradient $\nabla \phi_i$ requires only $\frac{1}{N}$ backward propagations. Hence, as in our implementation $I_{N_{k+1,g}}\subseteq  I_{\Nkpp}$, the computational cost of SIRTR at each iteration $k$ is determined by $\frac{N_{k+1}^t+N_{k+1,g}}{N}$ propagations. }

For all experiments in this section, we run SIRTR {with $x_0=(0,0,\ldots,0)^T$ as initial guess}, and stop it when either a maximum of $1000$ iterations is reached or a maximum of $500$ full function evaluations is performed or {  the condition 
\begin{equation}\label{eq:stop}
|f_{N_k}(x_k)-f_{N_{k-1}}(x_{k-1})|\leq \epsilon|f_{N_{k-1}}(x_{k-1})|+\epsilon,
\end{equation}
with $\epsilon=10^{-3}$, holds for a number of consecutive successful  iterations such that the computational effort is equal to the
effort needed in three iterations with full function and gradient evaluations}.

Since the selection of sets $I_{\Nkpp}$ and $I_{N_{k+1,g}}$  for computing $f_{\Nkpp}(x_k)$  and $g_k$  is
random, we perform $50$ runs of SIRTR for each test problem.
Results are reported in tables where the headings of the columns have the following meaning: {\tt cost} is the  overall number of full function and gradient evaluations averaged over the 50 runs, 
 {\tt err} is the classification error given in \eqref{class-err} averaged over the 50 runs,  {\tt sub}  the number of runs  where the method is stopped before reaching full accuracy in function evaluations.

In a first set of experiments, we investigate the choice of $N_{k+1,g}$
by varying the factor $c\in(0,1]$ in \eqref{eq:nkg}. In particular, letting  $\widetilde{c}=1.2$ in \eqref{eq:ntilde}, $\mu=100/N$ in \eqref{eq:nk1} and  $N_0=\lceil0.1N\rceil$ as in \cite{bkm}, we test the values $c\in\{0.1,0.2,1\}$. The results obtained are reported in Table \ref{tab:1}. We note that the classification error slightly varies with respect to the choice of $N_{k+1,g}$, and that selecting $N_{k+1,g}$ as a small fraction of $\Nkpp $  is quite convenient from  a computationally point of view. By contrast, the choice $N_{k+1,g}= \Nkpp $ leads to the largest computational costs without providing a significant gain in accuracy. Besides the cost per iteration, equal to 
$\frac{2\Nkpp}{N}$   in this latter case, we observe that full accuracy in function evaluations is reached   very often especially for certain datasets, see e.g., {\sc cina0}, {\sc cod-rna}, {\sc covertype},
{\sc ijcnn1}, {\sc phishing}, {\sc real-sim}.   Remarkably, the results in  Table \ref{tab:1} highlight that random models compare favourably with respect to cost and classification errors.  } 

\begin{table}[t!]
\begin{center}
\begin{tabular}{c|ccc|ccc|ccc}
$N_{k+1,g}$  &  \multicolumn{3}{|c|}{$\lceil 0.1 { \Nkpp} \rceil$} &  \multicolumn{3}{|c|}{$\lceil 0.2  { \Nkpp} \rceil$} & \multicolumn{3}{|c}{$ { \Nkpp}$} \\ \hline
 & {\tt cost} &  {\tt err} & {\tt sub} & {\tt cost}  & {\tt {err}} & {\tt sub} & {\tt cost} & {\tt {err}} & {\tt sub} \\\hline
{\sc a8a} 	 & 20 & 0.170 & 15 & 	19 & 0.171 & 19 & 	22 & 0.173 & 29 \\ 
{\sc a9a} 	 & 20 & 0.167 & 12 & 	17 & 0.169 & 18 & 	19 & 0.172 & 13 \\ 
{\sc cina0} 	 & 72 & 0.146 & 0 & 	84 & 0.140 & 0 & 	116 & 0.158 & 1 \\ 
{\sc cod-rna} 	 & 44 & 0.109 & 0 & 	42 & 0.106 & 1 & 	45 & 0.119 & 0 \\ 
{\sc covtype} 	 & 22 & 0.425 & 4 & 	19 & 0.424 & 8 & 	20 & 0.435 & 5 \\ 
{\sc htru2} 	 & 30 & 0.024 & 7 & 	25 & 0.024 & 13 & 	32 & 0.024 & 16 \\ 
{\sc ijcnn1} 	 & 22 & 0.087 & 0 & 	20 & 0.088 & 0 & 	20 & 0.086 & 0 \\ 
{\sc mnist2} 	 & 22 & 0.154 & 10 & 	25 & 0.151 & 12 & 	29 & 0.152 & 18 \\ 
{\sc phishing} 	 & 48 & 0.105 & 0 & 	43 & 0.108 & 0 & 	48 & 0.119 & 0 \\ 
{\sc real-sim} 	 & 56 & 0.268 & 0 & 	56 & 0.270 & 0 & 	57 & 0.294 & 0 \\ 
{\sc w7a} 	 & 15 & 0.079 & 22 & 	15 & 0.079 & 21 & 	16 & 0.079 & 34 \\ 
{\sc w8a} 	 & 13 & 0.080 & 25 & 	13 & 0.080 & 23 & 	17 & 0.080 & 28 
\end{tabular}
\end{center}
\caption{ Results with three different rules for computing the sample size $N_{k+1,g}$. 
}\label{tab:1}
\end{table}

Next, we show that SIRTR computational cost can be reduced by slowing down the growth rate of  $\Nkpp $. 
This task can be achieved controlling the growth of $\Ntilde$ which affects  $\Nkpp $  by means of  \eqref{eq:nk1}. 
Letting $c=0.1$, $\mu=100/N$ and  $N_0=\lceil0.1N\rceil$,
we consider the choices $\widetilde{c}\in\{1.05,1.1,1.2\}$ in \eqref{eq:ntilde}.
Average results are reported  in Table \ref{tab:2}. We can observe that {  the fastest  growth rate for $\Ntilde$ is generally more expensive than the other two choices, while  the classification error is similar for all the three choices}. Moreover,
significantly for $\widetilde c= 1.05$ {most} runs stopped before reaching full function accuracy. 

\begin{table}[t!]
\begin{center}
\begin{tabular}{c|ccc|ccc|ccc}
$\tilde{N}_{k+1}$  &  \multicolumn{3}{|c|}{$\min\{N,\lceil 1.05N_k\rceil\}$} &  \multicolumn{3}{|c|}{$\min\{N,\lceil 1.1N_k\rceil\}$} &  \multicolumn{3}{|c}{$\min\{N,\lceil 1.2N_k\rceil\}$} \\ \hline
 & {\tt cost} &  {\tt err} & {\tt sub} & {\tt cost}  & {\tt err} & {\tt sub} & {\tt cost} & {\tt err} & {\tt sub}\\\hline
{\sc a8a} 	 & 27 & 0.170 & 49 & 	18 & 0.170 & 44 & 	18 & 0.171 & 16 \\ 
{\sc a9a} 	 & 27 & 0.164 & 49 & 	18 & 0.164 & 38 & 	20 & 0.168 & 12 \\ 
{\sc cina0} 	 & 35 & 0.167 & 44 & 	44 & 0.163 & 13 & 	68 & 0.151 & 0 \\ 
{\sc cod-rna} 	 & 28 & 0.117 & 49 & 	38 & 0.108 & 17 & 	45 & 0.102 & 0 \\ 
{\sc covtype} 	 & 12 & 0.396 & 50 & 	13 & 0.392 & 48 & 	20 & 0.423 & 7 \\ 
{\sc htru2} 	 & 30 & 0.022 & 46 & 	24 & 0.022 & 26 & 	25 & 0.024 & 11 \\ 
{\sc ijcnn1} 	 & 21 & 0.089 & 50 & 	16 & 0.086 & 49 & 	22 & 0.088 & 0 \\ 
{\sc mnist2} 	 & 19 & 0.144 & 50 & 	18 & 0.144 & 42 & 	23 & 0.152 & 12 \\ 
{\sc phishing} 	 & 28 & 0.117 & 50 & 	30 & 0.110 & 23 & 	46 & 0.103 & 0 \\ 
{\sc real-sim} 	 & 36 & 0.254 & 50 & 	65 & 0.272 & 0 & 	57 & 0.267 & 0 \\ 
{\sc w7a} 	 & 26 & 0.078 & 50 & 	18 & 0.078 & 46 & 	14 & 0.079 & 22 \\ 
{\sc w8a} 	 & 20 & 0.079 & 50 & 	14 & 0.080 & 46 & 	13 & 0.080 & 26 \\ 
\end{tabular}
\end{center}
\caption{Results with three different rules for computing the sample size $\Ntilde$.}
\label{tab:2}
\end{table}

We now analyze three different values, $N_0\in \{\lceil0.001N\rceil, \lceil0.01N\rceil,\lceil0.1N\rceil\}$, for the initial  sample size $N_0$. We apply SIRTR with $\tilde{c}=1.05$ in \eqref{eq:ntilde}, $\mu=100/N$ in \eqref{eq:nk1}, and $c=0.1$ in \eqref{eq:nkg}. Results are reported in Table \ref{tab:3}. We can see that, reducing  $N_0$, the number of full function/gradient evaluations {can further reduce in some datasets}, and that for $N_0= \lceil0.01N\rceil$ the average classification error compares well with the error when $N_0= \lceil0.1N\rceil$; for instance, the best results {for most datasets} are obtained by shrinking $N_0$ to $1\%$ of the maximum sample size. {We conclude pointing out that  most of the runs are performed without reaching full precision in function evaluation.} 

{As a further confirmation of the efficiency of SIRTR, in Table \ref{tab:4} we report the sample sizes obtained on average at the stopping iteration of SIRTR with parameters setting $N_0=\lceil 0.01 N\rceil$, $N_{k+1,g}=\lceil 0.1 { \Nkpp}\rceil$, 
$\widetilde{N}_{k+1}=\min\{N,\lceil 1.05 N_k\rceil \}$, $\mu=100/N$. More specifically, for each dataset, we show the mean value 
$\overline{N}_{{\rm fin}}$ obtained by averaging the sample sizes $N_{{\rm fin},i}$, $ 1\le i\le 50$, used at the final iteration of SIRTR, the { relative standard deviation $s=\frac{1}{\overline{N}_{{\rm fin}}} \sqrt{\frac{\sum_{i=1}^{50}(N_{{\rm fin},i}-\overline{N}_{{\rm fin}})^2}{50}}$} as a measure of dispersion of the final sample sizes with respect to the mean value, and the minimum and maximum sample sizes $N_{{\rm fin}}^{\min},N_{{\rm fin}}^{\max}$ observed at the final iteration out of the 50 runs. From the reported values, we deduce that SIRTR terminates with a final sample size which is much smaller, on average, than the maximum sample size $N$. 

\begin{table}[t!]
\begin{center}
\begin{tabular}{c|ccc|ccc|ccc}
$N_0$  &  \multicolumn{3}{|c|}{$\lceil 0.001N\rceil$} &  \multicolumn{3}{|c|}{$\lceil 0.01N\rceil$} &  \multicolumn{3}{|c}{$\lceil 0.1N\rceil$} \\ \hline
 & {\tt cost} &  {\tt err} &{\tt sub}  & {\tt cost}  & {\tt err} & {\tt sub} & {\tt cost} & {\tt {err}} & {\tt sub}\\\hline
{\sc a8a} 	 & 30 & 0.182 & 50 & 	30 & 0.169 & 47 & 	28 & 0.170 & 50 \\ 
{\sc a9a} 	 & 27 & 0.177 & 50 & 	28 & 0.165 & 50 & 	25 & 0.165 & 50 \\ 
{\sc cina0} 	 & 43 & 0.111 & 37 & 	33 & 0.133 & 43 & 	34 & 0.162 & 44 \\ 
{\sc cod-rna} 	 & 4 & 0.412 & 50 & 	25 & 0.194 & 50 & 	29 & 0.114 & 48 \\ 
{\sc covtype} 	 & 6 & 0.406 & 50 & 	8 & 0.403 & 50 & 	12 & 0.406 & 50 \\ 
{\sc htru2} 	 & 38 & 0.036 & 40 & 	35 & 0.021 & 43 & 	31 & 0.021 & 47 \\ 
{\sc ijcnn1} 	 & 24 & 0.095 & 50 & 	25 & 0.095 & 50 & 	19 & 0.091 & 50 \\ 
{\sc mnist2} 	 & 18 & 0.185 & 50 & 	20 & 0.160 & 50 & 	21 & 0.143 & 50 \\ 
{\sc phishing} 	 & 4 & 0.410 & 50 & 	28 & 0.163 & 48 & 	29 & 0.118 & 50 \\ 
{\sc real-sim} 	 & 4 & 0.188 & 50 & 	5 & 0.166 & 50 & 	35 & 0.254 & 50 \\ 
{\sc w7a} 	 & 28 & 0.077 & 50 & 	27 & 0.077 & 50 & 	25 & 0.078 & 50 \\ 
{\sc w8a} 	 & 23 & 0.078 & 50 & 	23 & 0.079 & 50 & 	20 & 0.079 & 50 
\end{tabular}
\end{center}
\caption{Results with three different initial sample sizes $N_0$.}
\label{tab:3}
\end{table}

\begin{table}[t!]
{
\begin{center}
\begin{tabular}{c|c|c|c|c|c|c}
& $N$ & $N_0$ & $\bar{N}_{fin}$& $s$ & $N_{fin}^{\min}$ & $N_{fin}^{\max}$ \\
\hline
{\sc a8a} 	 & 15888 & 159 & 10353 & 0.17 & 7407 & 13309 \\ 
{\sc a9a} 	 & 22793 & 228 & 13637 & 0.22 & 6718 & 18730 \\ 
{\sc cina0} 	 & 10000 & 100 & 7603 & 0.23 & 4771 & 10000 \\ 
{\sc cod-rna} 	 & 7739 & 78 & 3210 & 0.74 & 578 & 7054 \\ 
{\sc covtype} 	 & 464810 & 4649 & 54762 & 0.32 & 33057 & 100341 \\ 
{\sc htru2} 	 & 10000 & 100 & 7902 & 0.22 & 3923 & 10000 \\ 
{\sc ijcnn1} 	 & 49990 & 500 & 26966 & 0.23 & 15408 & 43508 \\ 
{\sc mnist2} 	 & 60000 & 600 & 22928 & 0.34 & 4383 & 45684 \\ 
{\sc phishing} 	 & 7739 & 78 & 3926 & 0.63 & 578 & 7739 \\ 
{\sc real-sim} 	 & 50617 & 507 & 3721 & 0.034 & 3604 & 4174 \\ 
{\sc w7a} 	 & 17285 & 173 & 10334 & 0.23 & 5802 & 14674 \\ 
{\sc w8a} 	 & 34825 & 349 & 17244 & 0.19 & 9005 & 26360 
\end{tabular}
\end{center}

\caption{Average sample size $\overline{N}_{{\rm fin}}$ obtained at the final iteration, relative standard deviation $s$, minimum and maximum sample sizes $N_{{\rm fin}}^{\min},N_{{\rm fin}}^{\max}$ observed at the final iteration. Parameters setting: $N_0=\lceil 0.01 N\rceil$, $N_{k+1,g}=\lceil 0.1 { \Nkpp}\rceil$, $\widetilde{N}_{k+1}=\min\{N,\lceil 1.05 N_k\rceil \}$, $\mu=100/N$.}
\label{tab:4}}
\end{table}

Finally, in Figures \ref{fig:a9a_card}-\ref{fig:mnist2_card}, we report the plots of the sample sizes  $\Nkpp$ and $\widetilde{N}_{k+1}$ with respect to the number of iterations, obtained by running SIRTR on the {\sc a9a} and {\sc mnist} datasets, respectively. In particular, we let either $\mu=100/N$ or $\mu=1$ in the update rule \eqref{eq:nk1}, $\tilde{c}=1.05$ in \eqref{eq:ntilde}, $c=0.1$ in \eqref{eq:nkg} and $N_0=\lceil 0.1 N\rceil$. Note that a larger $\mu$ allows for the decreasing of both  $\Nkpp$ and $\widetilde{N}_{k+1}$ in the first iterations, whereas a linear growth rate is imposed only in later iterations. This behaviour is due to the update condition \eqref{eq:nk1}, which naturally forces $\Nkpp$ to coincide with $\widetilde{N}_{k+1}$ when $\delta_k$ is sufficiently small. For both choices of $\mu$, we see that $\Nkpp$  can grow slower than $\widetilde{N}_{k+1}$ at some iterations, thus reducing the computational cost per iteration of SIRTR. 

\begin{figure}[t!]
\begin{tabular}{c@{}c}
\includegraphics[scale = 0.42]{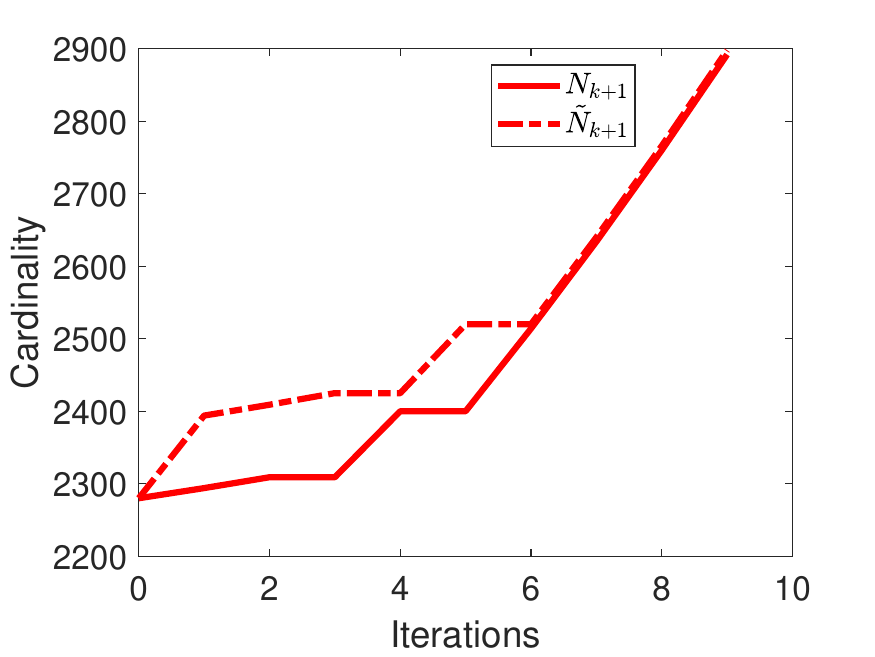} & \includegraphics[scale = 0.42]{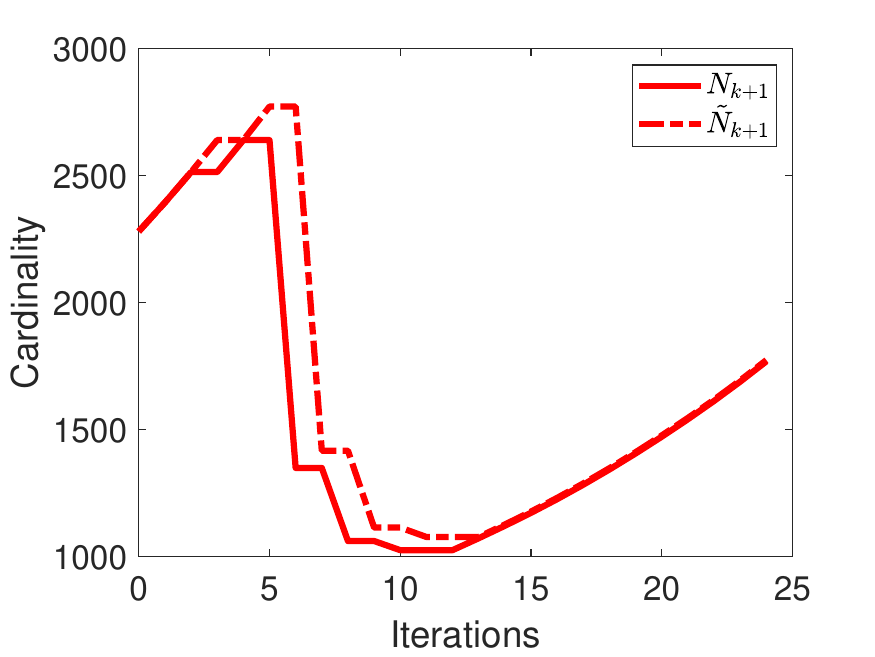} 
\end{tabular}
\caption{Dataset {\sc a9a}. Samples sizes $ N_{k+1} $ and $\widetilde{N}_{k+1}$ versus iterations with $\mu=100/N$ (left) and $\mu=1$ (right), respectively, obtained with a single run of SIRTR. Classification errors: {\tt err} = 0.187 with $\mu=100/N$, {\tt err} = 0.174 with $\mu=1$.}
\label{fig:a9a_card}
\end{figure}

\begin{figure}[t!]
\begin{tabular}{c@{}c}
\includegraphics[scale = 0.42]{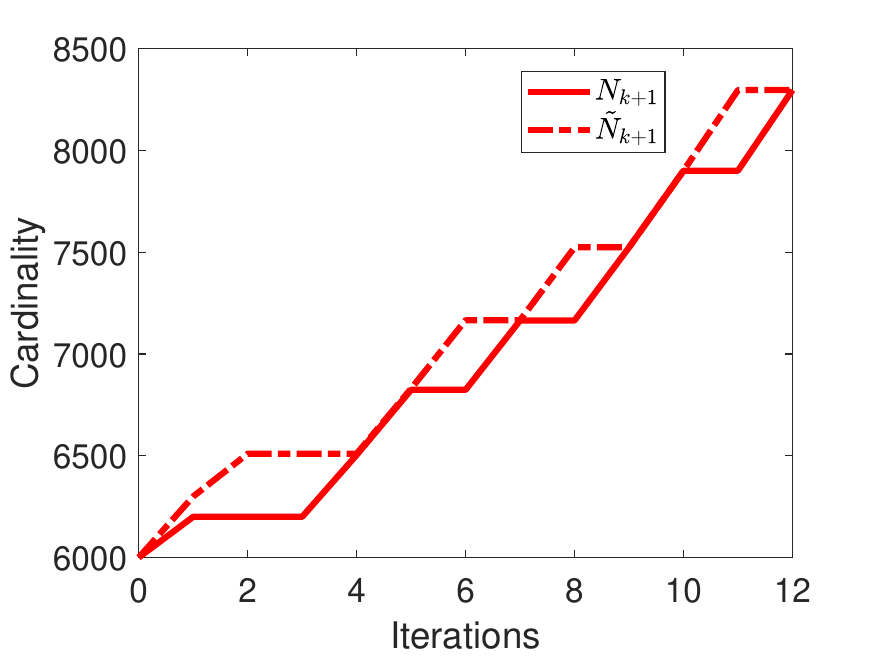} & \includegraphics[scale = 0.42]{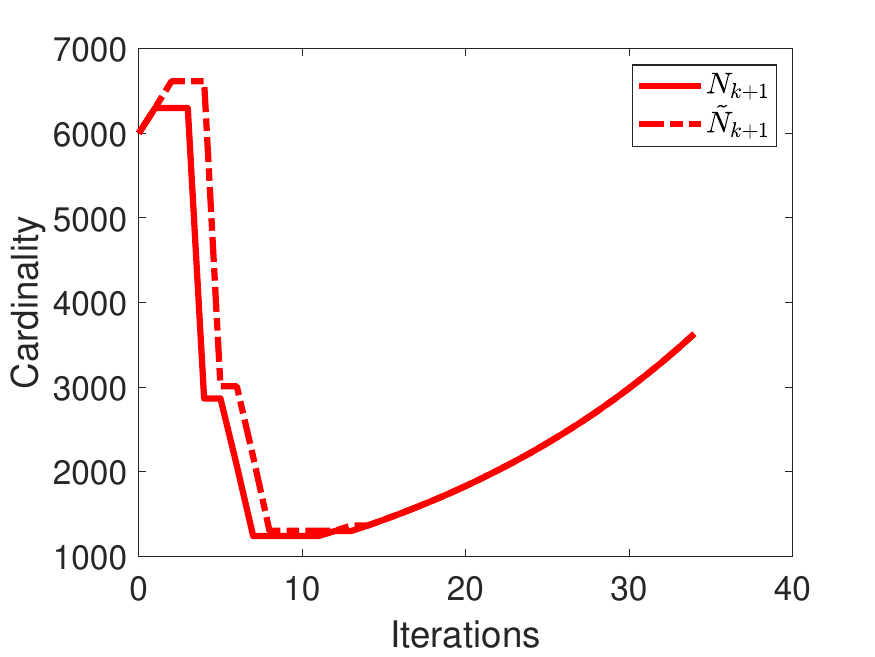} 
\end{tabular}
\caption{Dataset {\sc mnist}. Samples sizes $N_{k+1}$ and $\widetilde{N}_{k+1}$ versus iterations with $\mu=100/N$ (left) and $\mu=1$ (right), respectively, obtained with a single run of SIRTR. Classification errors: {\tt err} = 0.154 with $\mu=100/N$, {\tt err} = 0.167 with $\mu=1$.}
\label{fig:mnist2_card}
\end{figure}

\subsection{Comparison with TRish}

In this section we compare the performance of SIRTR with the so-called Trust-Region-ish algorithm (TRish) recently proposed in \cite{Curtis2019}. TRish is a stochastic gradient method based on a trust-region methodology. Normalized steps are used in a dynamic manner whenever the norm of the stochastic gradient is within a prefixed interval. In particular, the $k-$th iteration of TRish is given by 
\begin{equation*}\label{eq:TRish}
x_{k+1}=x_k-\begin{cases}
\gamma_{1,k}\alpha_kg_k, \quad & \text{if } \|g_k\|\in\left[0,\frac{1}{\gamma_{1,k}}\right)\\
\alpha_k\frac{g_k}{\|g_k\|}, \quad & \text{if }\|g_k\|\in\left[\frac{1}{\gamma_{1,k}},\frac{1}{\gamma_{2,k}} \right]\\
\gamma_{2,k}\alpha_kg_k, \quad &\text{if }\|g_k\|\in\left(\frac{1}{\gamma_{2,k}},\infty\right)
\end{cases}
\end{equation*}
where $\alpha_k>0$ is the steplength parameter, $0<\gamma_{2,k}<\gamma_{1,k}$ are positive constants, and $g_k\in\mathbb{R}^{n}$ is a stochastic gradient estimate. This algorithm has proven to be particularly effective on binary classification and neural network training, especially if compared with the standard stochastic gradient  algorithm \cite[Section 4]{Curtis2019}.

For our numerical tests, we implement TRish with subsampled gradients $g_k=\nabla f_{S}(x_k)$ defined in (\ref{subsample}). The steplength is constant, $\alpha_k=\alpha$, $\forall k\ge 0$, and $\alpha$ is chosen in the set  $ \{10^{-3},10^{-1},\sqrt{10^{-1}},1,\sqrt{10}\}$. Following the procedure in \cite[Section 4]{Curtis2019}, we use constant parameters $\gamma_{1,k}\equiv \gamma_1$, $\gamma_{2,k}\equiv \gamma_2$
and select $\gamma_1,\, \gamma_2$ as follows. First, Stochastic Gradient algorithm \cite{RM} is run with constant steplength equal to $1$;
second, the average norm $G$ of stochastic gradient estimates throughout the runs is computed; third  $\gamma_1,\, \gamma_2$ are set as $\gamma_1=\frac{4}{G}$, $\gamma_2=\frac{1}{2G}$.

%
%

\begin{figure}[t!]
\begin{center}
\begin{tabular}{c@{}c@{}c@{}}
\includegraphics[scale = 0.3]{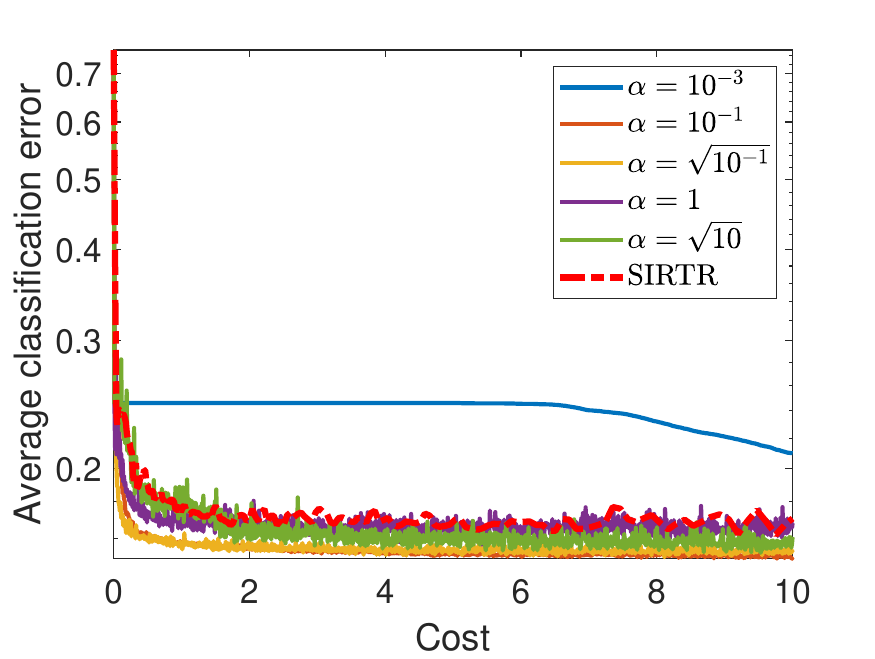} & \includegraphics[scale = 0.3]{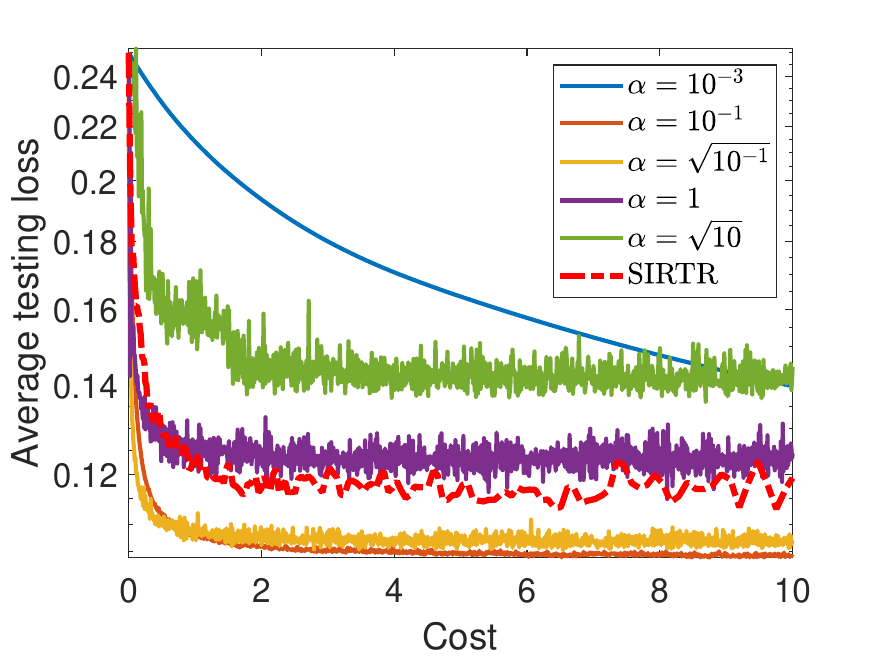} & \includegraphics[scale = 0.3]{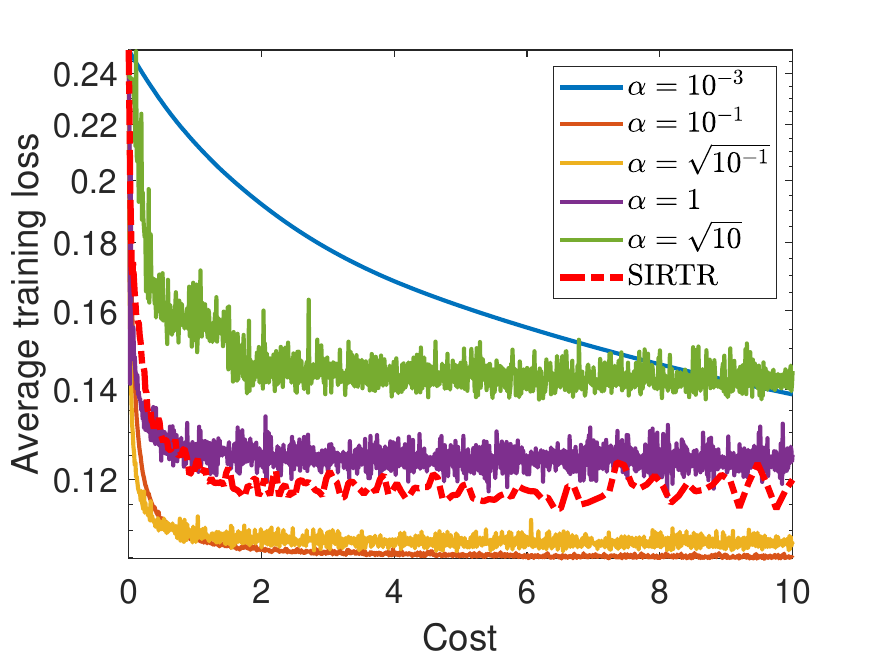}\\
\includegraphics[scale = 0.3]{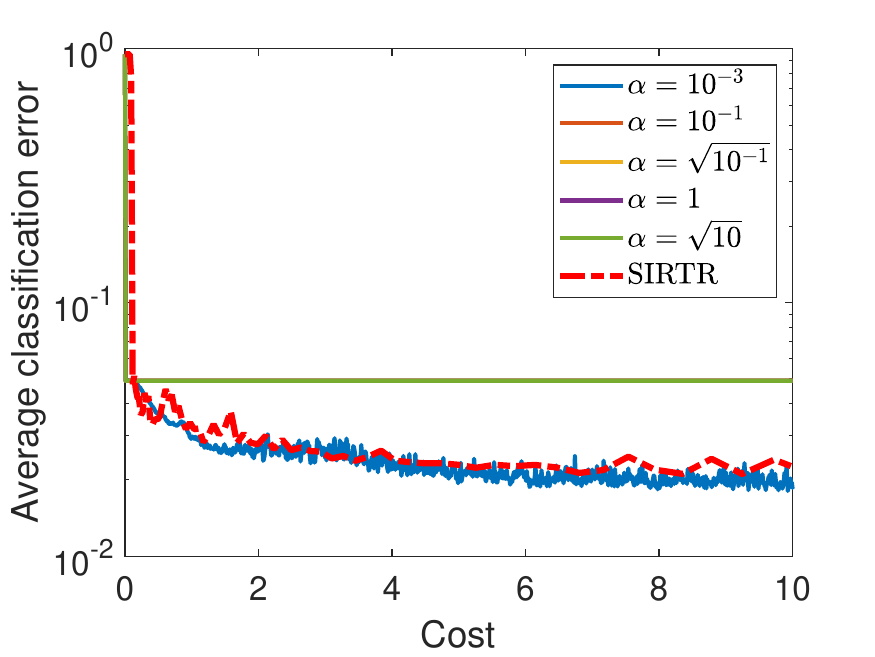} & \includegraphics[scale = 0.3]{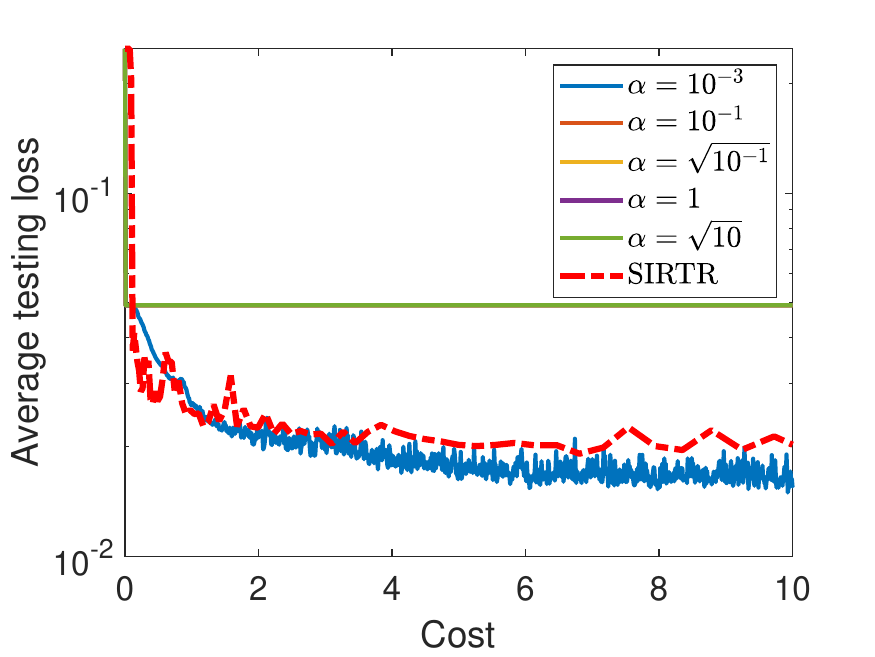} & \includegraphics[scale = 0.3]{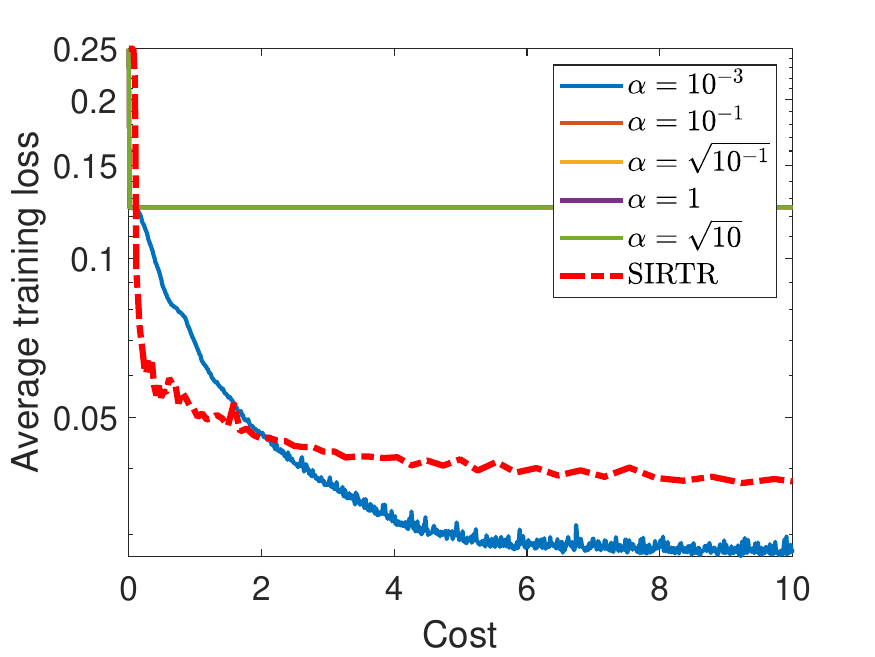}\\
\includegraphics[scale = 0.3]{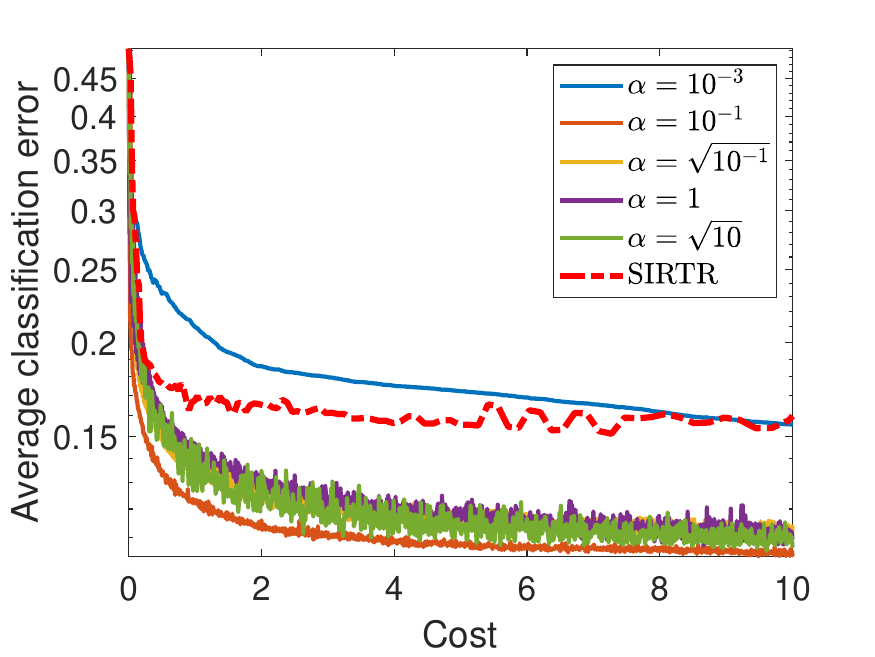} & \includegraphics[scale = 0.3]{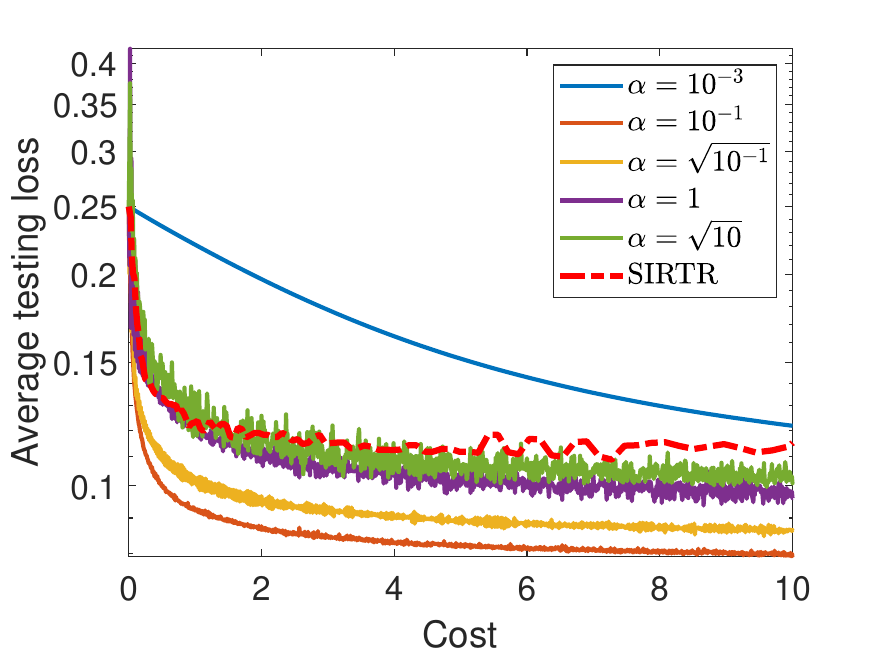} & \includegraphics[scale = 0.3]{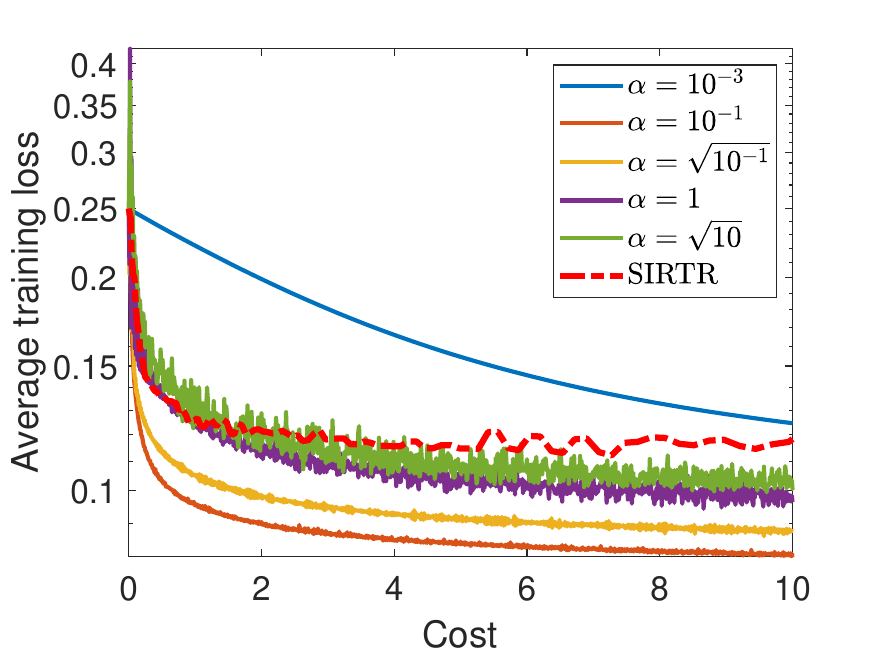}\\
\includegraphics[scale = 0.3]{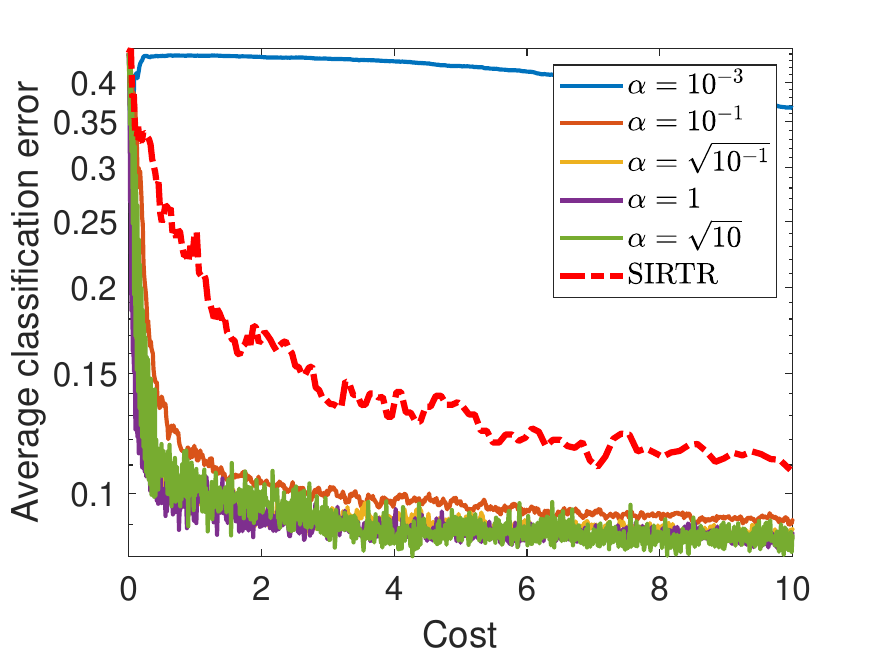} & \includegraphics[scale = 0.3]{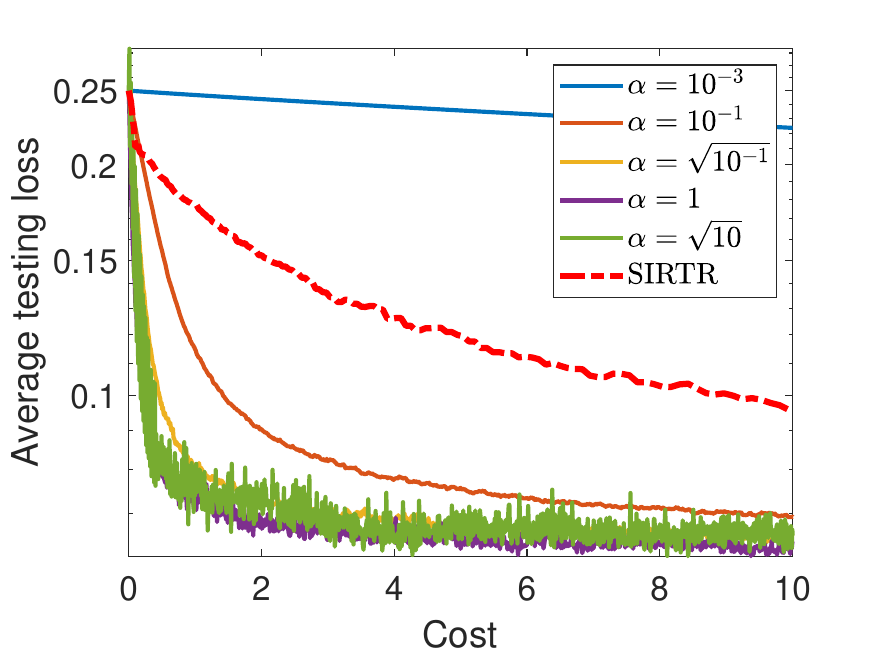} & \includegraphics[scale = 0.3]{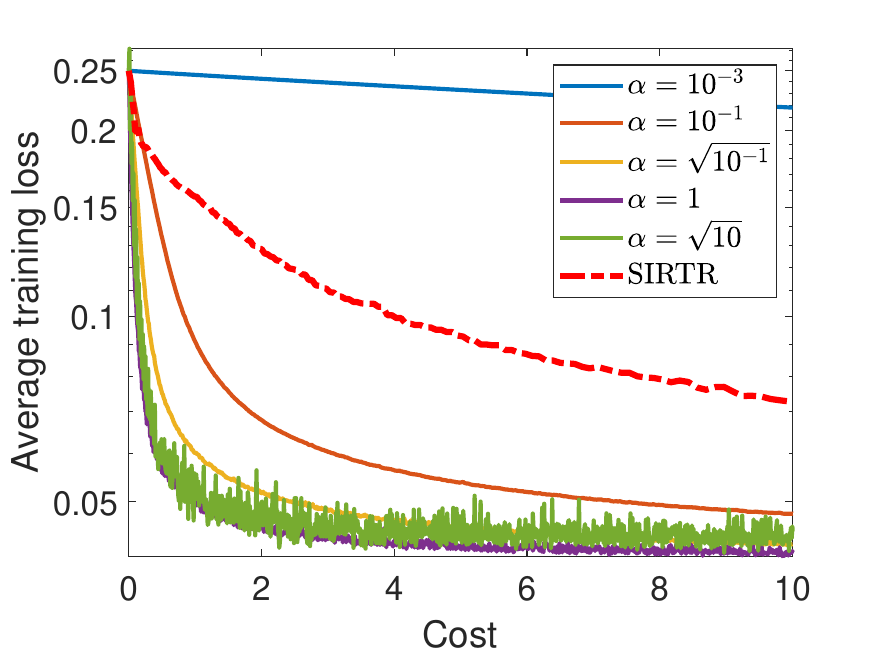}\\
\end{tabular}
\end{center}


\caption{ From top to bottom row: datasets {\sc a9a}, {\sc htru2}, {\sc mnist}, {\sc phishing}. From left to right: Average classification error, testing loss, and training loss versus epochs.}\label{fig:3}
\end{figure}

{First,} we compare  TRish  with SIRTR {on the nonconvex optimization problem \eqref{eq:sigmoid}, using {\sc a9a}, {\sc htru2}, {\sc mnist}, and {\sc phishing} as datasets (see Table \ref{test}).} Based on the previous section, {we equip SIRTR with} $N_0=\lceil 0.01 N\rceil$, $N_{k+1,g}=\lceil 0.1 { \Nkpp}\rceil$, $\widetilde{N}_{k+1}=\min\{N,\lceil 1.05 N_k\rceil \}$, $\mu=100/N$. In TRish, the sample size $S$ of the stochastic gradient estimates is 
$\lceil 0.01N\rceil $, which corresponds to the first sample size used in SIRTR. We run each algorithm for ten epochs on the datasets {\sc a9a} and {\sc htru2} using the null  initial guess. We perform   $10$ runs to report results on average. 

After tuning, the parameter setting for TRish was $ \gamma_1 \approx 34.5805$, $ \gamma_2 \approx 4.3226$ for {\sc a9a}, $ \gamma_1 \approx 57.9622$, $ \gamma_2 \approx 7.2453$ for {\sc htru2}, $ \gamma_1 \approx 23.4376$, $ \gamma_2 \approx 2.9297$ for {\sc mnist}, and $\gamma_1 \approx 50.6409$, $ \gamma_2 \approx 6.3301$ for {\sc phishing}. In Figure \ref{fig:3}, we report the decrease of the (average) classification error, training loss $f_N$ and testing loss,
$f_{N_T}(x)= \frac{1}{N_T}\sum_{i\in I_{N_T}} \phi_i(x)$, over the (average) number of full function and gradient evaluations required by the algorithms. From these plots, we can see that SIRTR performs comparably to the best implementations of TRish on {{\sc a9a}, {\sc htru2}, {\sc mnist}, while showing a good, though not optimal, performance on {\sc phishing}.} 

 In accordance to the experience in \cite{Curtis2019}, all  
parameters $\gamma_1$ and $\gamma_2$ and $\alpha$ are problem-dependent. { For instance}, the best performance of TRish is obtained with $\alpha=10^{-1}$ for {\sc a9a} and with $\alpha=10^{-3}$ for {\sc htru2}, respectively; by contrast, SIRTR performs { well with an unique setting of the parameters}  which is the key feature of adaptive stochastic 
optimization methods.

{As a second test, we compare the performance of SIRTR and TRish on a different nonconvex optimization problem arising from nonlinear regression. Letting $\{(a_i,b_i)\}_{i=1}^N$ denote the training set, where $a_i\in\IR^n$ and $b_i\in\IR$ represent the feature vector and the target variable of the $i $-th example, respectively, we aim at solving the following problem
\begin{equation}\label{eq:regression}
\min_{x\in\IR^n}f_N(x)=\frac 1 N\sum_{i=1}^N \left(b_i-h(a_i;x) \right)^2,
\end{equation}
where $h(\cdot;x):\IR^n\rightarrow \IR$ is a nonlinear prediction function. 

For this second test, we use the {\sc air} dataset \cite{UCI}, which contains $9358$ instances of (hourly averaged) concentrations of polluting gases, as well as temperatures and relative/absolute air humidity levels, recorded at each hour in the period March 2004 - February 2005 from a device located in a polluted area within an Italian city.

As in \cite{DeVito-et-al-08}, our goal is to predict the benzene (C6H6) concentration from the knowledge of $n=7$ features, including carbon monoxide (CO), nitrogen oxides (NO$_x$), ozone (O$_3$), non-metanic hydrocarbons (NMHC), nitrogen dioxide (NO$_2$), air temperature, and relative air humidity. First, we preprocess the dataset by removing examples for which the benzene concentration is missing, reducing the dataset dimension from $9357$ to $8991$. Then, we employ $70\%$ of the dataset for training ($N=6294$), and the remaining $30\%$ for testing ($N_T=2697$). Since the concentration values have been recorded hourly, this means that we use the data measured in the first $9$ months for the training phase, and the data related to the last $3$ months for the testing phase. Finally, denoting with $D=(d_{ij})\in \IR^{(N+N_T)\times n}$ the matrix containing all the dataset examples along its rows, and setting
$$
\begin{cases}
m_j = \min\limits_{i=1,\ldots,N+N_{T}} {d_{ij}}, \\
M_j = \max\limits_{i=1,\ldots,N+N_{T}} {d_{ij}}
\end{cases}, \quad j =1,\ldots,n,
$$
we scale all data values into the interval $[0,1]$ as follows
$$
d_{ij} = \frac{d_{ij}-m_{j}}{M_j-m_j}, \quad i=1,\ldots, N+N_T, \ j=1,\ldots,n.
$$
We apply SIRTR and TRish on problem \eqref{eq:regression}, where the prediction function $h(\cdot;x)$ is chosen as a feed-forward neural network based on a $7\times 5 \times 1$ architecture (see \cite{DeVito-et-al-08} and references therein), with the two hidden layers both equipped with the linear activation function, and the output layer with the sigmoid activation function. We equip the two algorithms with the same parameter  values employed  in the previous tests, and run them $10$ times for $10$ epochs, using a random initial guess in the interval $[-\frac{1}{2},\frac{1}{2}]$. 

In Figure \ref{fig:air_testing_training}, we report the decrease of the (average) training and testing losses provided by SIRTR and by TRish  with different choices of  the steplength $\alpha$, whereas in Figure \ref{fig:air_concentration} we show the benzene concentration estimations provided by the algorithms against  the true concentration. These results confirm
that the performances of SIRTR are comparable with those of TRish equipped with the best choice of the steplength  and show 
 the ability of SIRTR to automatically tune the  steplength  so as to obtain satisfactory results in terms of testing and training accuracy.

\begin{figure}[t!]
\begin{tabular}{c@{}c}
\includegraphics[scale = 0.42]{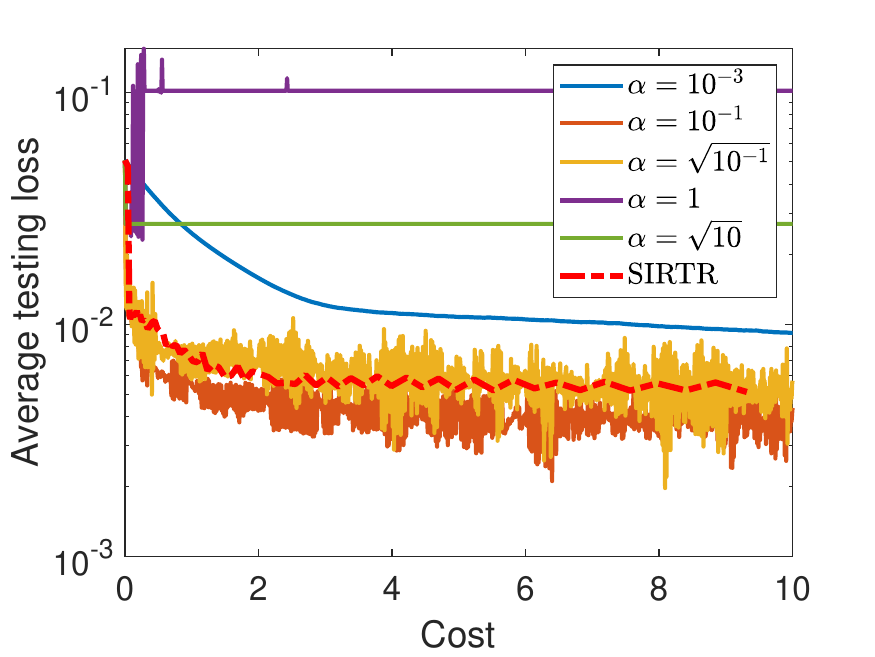} & \includegraphics[scale = 0.42]{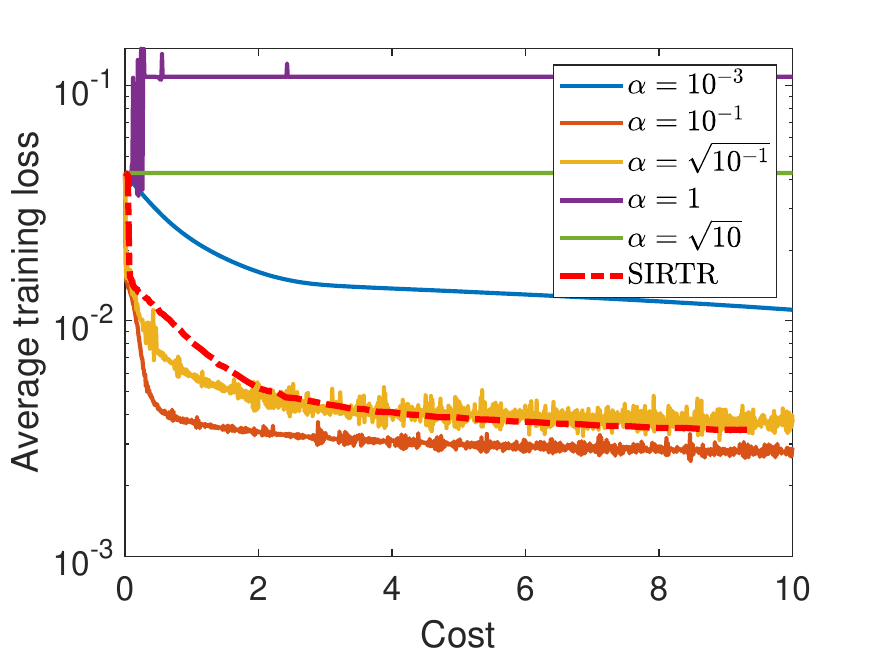} 
\end{tabular}
\caption{Dataset {\sc air}. Average testing loss (left) and training loss (right) versus epochs.}
\label{fig:air_testing_training}
\end{figure}

\begin{figure}[t!]
\begin{center}
\begin{tabular}{c}
\includegraphics[scale = 0.42]{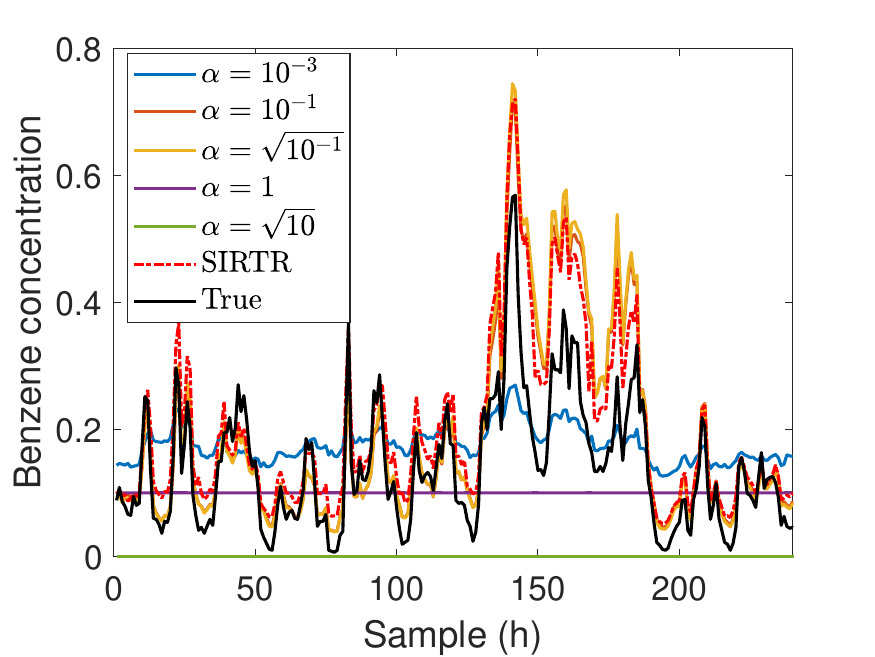}
\end{tabular}
\caption{Dataset {\sc air}. Estimated concentrations during $10$ days ($240$ hours) compared to the true concentration (black solid line).}
\label{fig:air_concentration}
\end{center}
\end{figure}

}

\section{Conclusions}
{ We proposed a stochastic gradient method coupled with a trust-region strategy and an inexact restoration approach for solving finite-sum minimization problems. Functions and gradients are subsampled and the batch size is governed by the inexact restoration approach and the trust-region acceptance rule. 
We showed the theoretical properties of the method and gave a worst-case complexity result on the expected number of iterations required to reach an approximate first-order optimality point.  Numerical experience shows that the proposed method provides good results  keeping the overall computational cost relatively low.}
\vskip 5pt\noindent
{\bf Data Availability}. The dataset CINA0 is no longer available in repositories but is available from the corresponding author on reasonable request. The other datasets analyzed during the current study are available in the repositories: 
\\
\url{http://www.csie.ntu.edu.tw/~cjlin/libsvm}, 
\\\url{ http://yann.lecun.com/exdb/mnist}, 
\\ \url{https://archive.ics.uci.edu/ml/index.php}
\vskip 5pt\noindent
{\bf Conflict of interest}. The authors have not conflict of interest to declare.

\end{document}